\documentclass{article}
\usepackage{amsmath,array,amsxtra,amssymb,amsbsy,graphicx}
\usepackage{supertabular,bm,color}
\usepackage{amsfonts}
\newcommand{\hgeom}[2]{\,{}_{#1\!}F_{#2}} % hypergeomertics
\def\ultra#1#2{P^{(#1)}_{#2}} % ultraspericals

\def\maxg{\mathsf r} %Minimum for |a|_p / \|g\|_p  in network

\def\NN{\mathbb{N}}
\def\CC{\mathbb{C}}

\def\cald{{\mathcal D}}

\def\calf{{\mathcal F}}
\def\calg{{\mathcal G}}
\def\calh{{\mathcal H}}

\def\calo{{\mathcal O}}
\def\calx{{\mathcal X}}

\def\cals{{\mathcal S}}

\def\K{{\mathcal K}}

\def\sfa{\mathsf {A}}
\def\sfb{\mathsf {B}}

\def\sfk{\mathsf {K}}
\def\sfl{\mathsf {L}}

\def\sfp{\mathsf {P}}
\def\sfq{\mathsf {Q}}
\def\seqb{{\mathsf b}}

\def\bfa{{\bf a}}
\def\bfx{{\bf x}}

\newcommand{\argmin}{\operatorname{argmin}}

\newcommand{\RR}{\mathbb{R}} \newcommand{\mq}{\mathrm{mq}}
 \newcommand{\sph}{\mathbb{S}}
\newcommand{\dist}{\operatorname{dist}} 
\newcommand{\spn}{\operatorname{span}}
\newcommand{\supp}{\operatorname{supp}}
\newtheorem{theorem}{Theorem}[section]
\newtheorem{lemma}[theorem]{Lemma}
\newtheorem{cor}[theorem]{Corollary}

\newtheorem{prop}[theorem]{Proposition}
\newtheorem{remark}[theorem]{Remark} \newenvironment{proof}{{\noindent
    \bf Proof:}}{\hfill$\Box$\bigskip}

\numberwithin{equation}{section}

\title{$L^p$ Bernstein Estimates and Approximation by Spherical Basis
  Functions
}

\author{H.~N.~Mhaskar\thanks{ Department of Mathematics, California
    State University, Los Angeles,CA 90032, USA. Research supported
    by grant DMS-0605209 from the National Science Foundation and
    grant W911NF-04-1-0339 from the U.S. Army Research Office.}, 
F. J. Narcowich\thanks{ Department of Mathematics, Texas A\&M
    University College Station, TX 77843, USA. Research
    supported by grants DMS-0504353 and DMS-0807033 from the National
    Science Foundation.}, 
J. Prestin\thanks{Institute of Mathematics,
    University of L\"ubeck, Wallstrasse 40, 23560, L\"ubeck, Germany},
J. D. Ward\thanks{ Department of Mathematics, Texas A\&M University
    College Station, TX 77843, USA. Research supported by
    grants DMS-0504353 and DMS-0807033 from the National Science
    Foundation.}  }

\date{}
\begin{document}
\maketitle

\begin{abstract}
  The purpose of this paper is to establish $L^p$ error estimates, a
  Bernstein inequality, and inverse theorems for approximation by a
  space comprising spherical basis functions located at scattered
  sites on the unit $n$-sphere. In particular, the Bernstein
  inequality estimates $L^p$ Bessel-potential Sobolev norms of
  functions in this space in terms of the minimal separation and the
  $L^p$ norm of the function itself. An important step in its proof
  involves measuring the $L^p$ stability of functions in the
  approximating space in terms of the $\ell^p$ norm of the
  coefficients involved. As an application of the Bernstein
  inequality, we derive inverse theorems for SBF approximation in the
  $L^P$ norm. Finally, we give a new characterization of Besov spaces
  on the $n$-sphere in terms of spaces of SBFs.
  \end{abstract}

\section{Introduction}

Various applications in meteorology, cosmology, and geophysics require
a modeling of functions based on \emph{scattered data} collected on
(or near) a sphere; i.e., when one does not have any control on where
the data sites are located
\cite{Freeden-etal-97-1,Freeden-etal-98-1,Freeden-etal-04-1}. On $\sph^n$, the unit sphere in $\RR^{n+1}$, $n\ge 1$, a popular method is to construct the required approximation from spaces of spherical basis functions (SBFs), which are kernels located at points in a discrete set  $X=\{\xi_j\}_{j=1}^N\in \sph^n$, the set of \emph{centers} or \emph{nodes}. 

A function $\phi : [-1,1]\to\RR$ is an SBF on $\sph^n$ if, in its expansion in ultraspherical polynomials $\ultra {\lambda_n} \ell$, $\lambda_n=\frac{n-1}{2}$, the Fourier-Legendre coefficients $\{\hat\phi(\ell)\}$ of $\phi$ are all positive; see section \ref{SBFs}
for details.  These $\phi$ are to be used as kernels of the form $\phi(x\cdot y)$, $x,y\in \sph^n$, $x\cdot y$ being the usual ``dot'' product. The approximation space here is the span 
\[
\calg_{\phi,X}:= \spn\{\phi(x\cdot\xi)\}_{\xi\in X}.
\]
Following usage common in the neural network community, we will say that a function $g\in \calg_{\phi,X}$ is an \emph{SBF network} associated with $\phi$. The SBF $\phi$ is sometimes called an \emph{activation function} or a \emph{neuron}, but we will not use these terms here.

Such $\phi$ may have singular behavior. This is the case for certain thin-plate splines; $(1-x\cdot y)^{-1/2}$ is an SBF in $\sph^n$, $n\ge 2$, for instance. However, when they are continuous, they are positive definite in Schoenberg's sense \cite{Schoenberg-42-1}. In that case the interpolation
matrix $[\phi(\xi_i\cdot\xi_j)]$ is positive definite, and it is
possible to use SBFs to interpolate data given at the points in $X$.

The focus of this paper is approximation. To handle noisy data, both least squares and quasi-interpolants have been used for many years. More recently, the issue in many meshless numerical methods for solving PDEs is how well a network approximates a solution to the PDE. Singular SBFs should prove useful in probing for a corresponding singularity in solutions. 

To be effective, though, such methods require knowing the degree of approximation in various spaces, especially the $L^p$, $1\le p\le \infty$. The $L^2$ case for SBFs $\phi$ with $\hat\phi(\ell) \sim (\ell+1)^{-\beta}$, $\beta>n/2$ has recently been investigated in \cite{Narcowich-etal-07-1}, with nearly optimal rates being attained by interpolatory networks. The known estimates on the degree of approximation in the case of $L^p$, $p\ne 2$ provided by interpolatory networks are not asymptotically optimal. This has lead to the development of other approximation tools  \cite{Mhaskar-etal-99-1, Mhaskar-06-1,Narcowich-etal-06-1}, involving SBFs or spherical harmonics, in $L^p$, $1\le p\le \infty$. A central step in obtaining approximation rates in $L^2$ was establishing a Bernstein estimate, which was the used to get an inverse approximation theorem.  

The paper has three main goals. The first is to derive an $L^p$  Bernstein
inequality, for $1\le p\le \infty$; namely, $\|g\|_{H^p_\gamma} \le Cq^{-\gamma}\|g\|_p$, $0<\gamma<c_\phi$. Here $H^p_\gamma$ is a Bessel-potential Sobolev space \cite{Strichartz-83-1,Triebel-86-1}; it measures derivatives of $g$ (cf. section~\ref{bessel-sobolev}). The quantity $q$ is a half the minimal separation of points in $X$; $q^{-1}$ plays the role of a Nyquist frequency. 

The second is to obtain is to obtain $L^p$ error estimates, $1\le p\le \infty$, for approximating a function by networks in $\calg_{\phi, X}$. We combine these direct (Favard-Jackson) estimates with the Bernstein inequalities to provide new characterizations of Besov spaces on $\sph^n$, characterizations that use rates of approximation from the $\calg_{\phi,X}$. The Bernstein estimates are then used to establish inverse theorems and obtain nearly optimal rates of approximation.

The third is to show that the results gotten here will apply for nearly all of the SBFs of interest. In particular,  they apply to various RBFs restricted to the sphere -- the thin-plate splines and Wendland functions, whose Fourier-Legendre coefficients have algebraic decay, and also Gaussians and multiquadrics, whose coefficients decay faster than algebraically. SBFs in the latter class are well known to be difficult to treat.

The paper is organized this way. Section~\ref{background} reviews various geometric quantities, such as the set of centers, mesh norm, and so on. It also discusses spherical harmonics and the Bessel-potential Sobolev spaces. Section~\ref{SBFs} discusses SBFs, their Fourier-Legendre expansions, and deals in detail with the SBFs mentioned earlier, along with ones corresponding to certain Green's functions that play a significant role in the paper. It is here that we will show that nearly all of the SBFs of interest have the properties necessary for our results will hold. We also mention that we obtain precise asymptotic expressions for the Fourier-Legendre coefficients in the case of the Wendland functions. 

The strategy for establishing the Bernstein inequality, which will be detailed below, consists of two key components:  $L^p$ approximation results for functions in $\calg_{\phi,X}$ by means of spherical polynomials, and $L^p$ stability estimates; these are developed in sections~\ref{approximation}~and~\ref{stability}, respectively. The approximation results are based on Marcinkiewicz-Zygmund inequalities developed in \cite{Mhaskar-etal-01-1,Mhaskar-etal-01-2, Narcowich-etal-06-1}, as well as frame results from \cite{Narcowich-etal-06-1}. The stability results, which are of interest in their own right, are for all $L^p$ -- not just for interpolation with continuous SBFs. To obtain them, we introduce a \emph{stability ratio}, which provides some measure of the extent to which a finite set in $L^p$ is linearly independent. 

In section~\ref{bernstein_inverse_theorems}, the results of the previous two sections are combined to yield $L^p$ Bernstein inequalities (section~\ref{bernstein_inequalities}), direct theorems for approximation by networks in $\calg_{\phi,X}$ (section~\ref{direct_thms}), characterizations of Besov spaces on $\sph^n$ (section~\ref{besov_spaces}),  and inverse theorems for $L^p$ functions approximated at given rates by SBF networks (section~\ref{inverse_theorems}).

\paragraph{Strategy} Let $g$ be an SBF network in $\calg_{\phi,X}
\subset H^p_\gamma(\sph^n)$, so that it has the form
\[
g(\bfx)=\sum_{\xi\in X}a_\xi \phi(\bfx \cdot \xi).
\]
One of our main goals is to obtain an $L^p$ Bernstein inequality for such
networks; that is,  a bound of the form $\|g\|_{H^p_\gamma} \le Cq^{-\gamma}\|g\|_p$, where the norms are those appropriate for $\sph^n$ and $\gamma>0$ is bounded above by a
constant depending on $\phi$ and $p$.

Our strategy involves approximating $g$ by degree $L$ spherical polynomials on
$\sph^n$, where $L\sim q^{-1}$. Now, for fixed $L$ and any $S$, there is a Bernstein
inequality, $\|S\|_{H^p_\gamma}\le CL^\gamma \|S\|_p$, which is found in Theorem~\ref{bernstein-nikolskii-poly-ineqs}. Using it and manipulations involving the triangle inequality,
one has that
\[
\|g\|_{H^p_\gamma} \le \|S\|_{H^p_\gamma}+\|g-S\|_{H^p_\gamma}\le
CL^\gamma \|S\|_p+\|g-S\|_{H^p_\gamma},
\]
which holds for given $L$ and any $S$. 

Obtaining an appropriate polynomial $S$ is crucial to the argument. To
do that, we will use the frame operators introduced in
\cite{Narcowich-etal-06-1} and discussed in more detail in section
\ref{frames} below. In particular, we need reconstruction operators
$\sfb_J$, with $J\sim \log_2 L$.  These rotationally invariant
operators have other very useful approximation properties, which are
given in Proposition \ref{B_J_approx_thm}. They take $L^p$ spaces and
the space of continuous function boundedly into spherical polynomials
having degree $\calo(2^J)$. Consequently, with $S=\sfb_J g$, we have
$\|S\|_p\le C\|g\|_p$, and also
\[
\|g\|_{H^p_\gamma} \le C2^{\gamma J} \|g\|_p+\|g-\sfb_J
g\|_{H^p_\gamma} = C2^{\gamma J} \|g\|_p + \frac{|a|_p}{\|g\|_p} \cdot
\frac{\|g-\sfb_J g\|_{H^p_\gamma}}{|a|_p} \cdot \|g\|_p
\]
where $|a|_p = \left(\sum_{\xi \in X} |a_\xi|^p \right)^{1/p}$ is the
$p$-norm of $a=\{a_\xi\}_{\xi\in X}$.

The functions $\{ \phi((\cdot)\cdot \xi)\}_{\xi\ni X}$ are linearly
independent and form a basis for $\calg$, and so the pairing
$a\leftrightarrow g$ is bijective.  Since $\calg$ has finite dimension
$| \calg |$, the ratio
\begin{equation}
\label{s-ratio}
\maxg_{\calg,\, p}:= \max_{\calg \ni g \neq 0}\frac{|a|_p}{\|g\|_p}
\end{equation}
is finite; it will be called the $p$-norm \emph{stability ratio} of
the network $\calg=\calg_{\phi,X}$. This ratio is similar to a
condition number in interpolation, but for $L^p$. With it, the
inequality directly above becomes
\begin{equation}
\label{basic_ineq}
\|g\|_{H^p_\gamma} \le \left( C2^{\gamma J} + C'\maxg_{\calg,p} 
\left(\frac{\|(I-\sfb_J)g\|_{H^p_\gamma}}{|a|_p} \right) \right)\|g\|_p.
\end{equation}

To obtain the desired Bernstein inequality, we require two bounds: the
first on $\|(I-\sfb_J)g\|_{H^p_\gamma}/|a|_p$ and the second on
$\maxg_{\calg,p}$. The first bound relies only on approximation
results; these we cover in section \ref{approximation}. The second is
a bound on the stability ratio. This bound requires a more detailed
analysis involving both the geometry of $X$ and properties of
$\phi$. It is carried out in section \ref{stability}.

An interesting point is the that the two bounds make different demands
on the properties required for $\phi$. This makes the analysis of both
bounds subtle. Fortunately, the common demands are satisfied by large
classes of SBBs, including restrictions to $\sph^n$ of the most common
RBFs -- the thin-plate splines, Wendland functions, Gaussians, Hardy
multiquadrics, and others.

\section{Background}\label{background}
\subsection{Background and notation for $\sph^n$}
\label{notation}

\paragraph{Centers and decompositions of $\sph^n$.} 
Let $X$ be a finite set of distinct points in $\sph^n$; we will call
these the \emph{centers}. For $X$, we define these quantities: \emph{mesh norm}, $h_X=\sup_{y\in
  \sph^n} \inf_{\xi\in X} d(\xi,y)$, where $d(\cdot,\cdot)$ is the
geodesic distance between points on the sphere; the \emph{separation
  radius}, $q_X=\frac12 \min_{\xi\neq\xi'}\, d(\xi,\xi')\,$; and the
\emph{mesh ratio}, $\rho_X:=h_X/q_X\ge 1$.

For $\rho\ge 1$, define $\calf_\rho=\calf_\rho(\sph^n)$ be the family
of all sets of centers $X$ with $\rho_X\le \rho\,$. We say that $X$ is
\emph{$\rho$-uniform} if $X\in \calf_\rho$.  For every $\rho\ge 2$,
$\calf_\rho(\sph^n)$ is not only non empty, but it contains
\emph{nested} sequences of sets of centers for which $h_X$ becomes
arbitrarily small; precisely, the result is this:

\begin{prop}[{\cite[Proposition 2.1]{Narcowich-etal-07-1}}]\label{nesting_Xs}
  Let $\rho\ge 2$ and let $\calf_\rho$ be the corresponding
  $\rho$-uniform family. Then, there exists a sequence of sets $X_k\in
  \calf_\rho$, $k=0,1,\ldots$, such that the sequence is nested,
  $X_k\subset X_{k+1}$, and such that at each step the mesh norms
  satisfy $\frac14 h_{X_k} < h_{X_{k+1}}\le \frac12 h_{X_k}$.
\end{prop}

We will need to consider a decomposition of $\sph^n$ into a finite
number of non-overlapping, connected regions $R_\xi$, each containing
an interior point $\xi$ that will serve for function evaluations as
well as labeling.  For example, if $\calx$ is the Voronoi tessellation
for a set of centers $X$, then we may take $R_\xi$ to be the region
associated with $\xi\in X$. In any case, we will let $X$ be the set of
the $\xi$'s used for labels and $\calx=\{R_\xi\subset \sph^n\,|\,
\xi\in X\}$. In addition, let $\|\calx\|=\max_{\xi\in X}
\{\mbox{diam}(R_\xi)\}$. 

\subsection{Spherical harmonics}
Let $n\ge 2$. Let $d\mu$ be the standard measure on the $n$-sphere,
and let the spaces $L^p(\sph^n)$, $1\le p\le \infty$, have their usual
meanings. In addition, let $\Delta_{\sph^n}$ denote the
Laplace-Beltrami operator on $\sph^n$. The eigenvalues of
$\Delta_{\sph^n}$ are $-\ell(\ell +n-1)$, $\ell=0,1,\ldots$. For $n\ge
2$ and $\ell$ fixed, the dimension of the eigenspace is
\begin{equation}
\label{H_dimension}
d_\ell^n =\frac{\ell+\lambda_n}{\lambda_n} \binom{\ell+n-2}{\ell}  
\ \stackrel{\ell\to\infty}{\sim}\  \frac{\ell^{n-1}}{\lambda_n
  (n-2)!},\ 
\lambda_n:=\frac{n-1}{2}.
\end{equation}
For $n=1$, the case of the circle, $d_0^1=1$ and $d_\ell^1=2$, $\ell\ge1$.

A spherical harmonic $Y_{\ell,m}$ is an eigenfunction of
$\Delta_{\sph^n}$ corresponding to the eigenvalue $-\ell(\ell +n-1)$
\cite{Mueller-66-1, Stein-Weiss-71-1}, where $m=1\ldots d_\ell^n$. The
set $\{Y_{\ell,m}\;\colon \ell=0,1, \ldots, m=1\ldots d_\ell^n\}$ is
orthonormal in $L^2(\sph^n)$. Denote by $\calh_\ell$ the span of the
spherical harmonics with fixed order $\ell$, and let
$\Pi_L=\bigoplus_{\ell=0}^L \calh_\ell$ be the span of all spherical
harmonics of order at most $L$.  The orthogonal projection $\sfp_\ell$
onto $\calh_\ell$ is given by
\begin{equation}
\label{proj_L}
\sfp_\ell f = \sum_{m=1}^{d_\ell^n} \langle f,Y_{\ell,m} \rangle
Y_{\ell,m}\,.
\end{equation}

We regard the sphere $\sph^n$ as being the unit sphere in $\RR^{n+1}$, and we let
the quantity $\xi\cdot \eta$ denote the usual ``dot'' product for
$\RR^{n+1}$. Using the addition formula for spherical harmonics, when $n\ge 2$, one can write the
kernel for this projection as
\begin{equation}
\label{kernel_proj_L}
P_\ell(\xi\cdot\eta) = \sum_{m=1}^{d_\ell^n} Y_{\ell,m}(\xi) \overline{
Y_{\ell,m}(\eta)} = \frac{\ell+\lambda_n}{\lambda_n\omega_n}\ultra
{\lambda_n} \ell (\xi\cdot \eta), \ \lambda_n := \frac{n-1}{2},
\end{equation}
where $\ultra {\lambda_n} \ell (\cdot)$ is the ultraspherical
polynomial of order $\lambda_n$ and degree $\ell$. Also, we have that
$\|\ultra {\lambda_n} \ell \|_\infty \le \ultra {\lambda_n} \ell
(1)=\frac{d_\ell^n \lambda_n}{\ell+\lambda_n}$. We will briefly
discuss these polynomials in section \ref{SBFs}, in connection
spherical basis functions. For $n=1$, $\lambda_1=0$. In that case, the
kernel for $\sfp_\ell$ has the form
\begin{equation}
\label{kernel_proj_L_n=1}
P_\ell(\xi\cdot\eta) = 
\left\{
\begin{array}{ll}
\frac{1}{2\pi}, &\ell=0\\[5pt]
\frac{1}{\pi}T_\ell(\xi\cdot \eta ), &\ell \ge 1,
\end{array}\right.
\end{equation}
where $T_\ell(\cdot)$ the degree-$\ell$ Chebyshev polynomial of the
first kind, which is a limiting case of the ultraspherical polynomials
\cite[Section 4.7]{Szego-75-1}.

We will also need to consider operators of the form
$\sum_{\ell=0}^\infty c_\ell \sfp_\ell$. The kernels for the
projections $\sfp_\ell$ then provide us with kernels
$\sum_{\ell=0}^\infty c_\ell P_\ell(\xi\cdot\eta)$, which may be
distributional.

\subsection{Bessel-potential Sobolev spaces} \label{bessel-sobolev}
The spherical harmonic $Y_{\ell,m}$ is an eigenfunction corresponding
to the eigenvalue $-\ell(\ell +n-1)=\lambda_n^2-(\ell+\lambda_n)^2$
for Laplace-Beltrami operator $\Delta_{\sph^n}$ on $\sph^n$. It
follows that $\ell+\lambda_n$ is an eigenvalue corresponding to the
eigenfunctions $Y_{\ell,m}\,, m=1\ldots d^n_\ell$, of the
pseudo-differential operator
\begin{equation}
\label{psiDO_Ln}
\sfl_n:=\sqrt{\lambda_n^2-\Delta_{\sph^n}} = 
\sum_{\ell=0}^\infty (\ell+\lambda_n)\sfp_\ell.
\end{equation}

Let $\gamma$ be real, $1\le p \le \infty$ and $n\ge 2$.  If $f$ is a distribution on $\sph^n$, define the Bessel-potential Sobolev spaces  $H^p_\gamma(\sph^n)$ \cite{Strichartz-83-1,Triebel-86-1} to be all $f$ such that
\begin{equation}\label{def-H}
\|f\|_{H^p_\gamma}:=
\Big\|\sum_{\ell=0}^\infty (\ell+\lambda_n)^{\gamma} \sfp_\ell f \Big\|_{L^p}
< \infty,
\end{equation}
where $\sfp_\ell$ is from (\ref{proj_L}). The notation we use here is
that of Triebel \cite{Triebel-86-1}. Strichartz \cite{Strichartz-83-1}
defined these spaces on an a complete Riemannian manifold, using the
equivalent operator $(1-\Delta_{\sph^n})^{\gamma/2}$ to do
so. One more thing: 

\begin{remark}\label{dom_L}
The space $H^2_\gamma(\sph^n)$ is the domain
of $\sfl_n^\gamma$ {\rm \cite[Theorem 4.4]{Strichartz-83-1}}, which implies
that $H^2_\gamma(\sph^n)$ is norm equivalent to the usual sobolev
space $W^\gamma_2(\sph^n)$. 
\end{remark}

\section{Spherical basis functions} \label{SBFs}

For any real $\lambda>0$, not just $\lambda_n=\frac{n-1}2$, the
ultraspherical polynomials satisfy the orthogonality relation,
\begin{equation}
\label{orthog_rel}
\int_{-1}^1\ultra {\lambda} \ell (x) \ultra {\lambda} k
(x)(1-x^2)^{\lambda- \frac12}dx = 
\frac{2^{1-\lambda}\pi \Gamma(\ell+2\lambda)}{(\ell+\lambda)\Gamma^2(
  \lambda) \Gamma(\ell+1)}  \delta_{k,\ell}.
\end{equation}
For the circle, we have $\lambda_1=0$. With $\ell\ge 1$, as $\lambda
\to 0$, the ratio $ \ultra {\lambda_n} \ell (\cdot)/\lambda$ converges
to $(2/\ell)T_\ell(\cdot)$, the degree-$\ell$ Chebyshev polynomial of
the first kind \cite[Section 4.7]{Szego-75-1}.

Consider a function $\phi$ in $L^p$ or $C$.  We will assume that
$\phi$ has the following expansion in the orthogonal set of
ultraspherical polynomials:
\begin{equation}\label{sbf_def}
\phi(\underbrace{\xi\cdot \eta}_{\cos\theta}) :=
\left\{
\begin{array}{ll}
\frac{1}{2\pi}\hat\phi(0) + 
\frac{1}{\pi}\sum_{\ell=1}^\infty
\hat\phi(\ell)\cos \ell \theta, &n=1,\\ [5pt]
\sum_{\ell=0}^\infty
\hat\phi(\ell)
\frac{\ell+\lambda_n}{\lambda_n\omega_n}\ultra
{\lambda_n} \ell (\cos \theta ), &n\ge 2.
\end{array}\right.
\end{equation}
where $\omega_n := \frac{2\pi^{\frac{n+1}{2}}}{\Gamma(\frac{n+1}{2})}$
is the volume of $\sph^n$.

Functions of this form are called \emph{zonal}. We will assume that
the series converges in at least a distributional sense. The
coefficients in the expansion are obtained via the orthogonality
relations in (\ref{orthog_rel}). These are given below.
\[
\frac{\ell+\lambda_n}{\lambda_n\omega_n}\hat\phi(\ell)=
\frac{(\ell+\lambda_n)\Gamma^2(\lambda_n)\Gamma(\ell+1)}
{2^{1-\lambda_n}\pi \Gamma(\ell+2\lambda_n)} \int_{-1}^1\phi(x)\ultra
{\lambda_n} \ell (x) (1-x^2)^{\lambda_n-\frac12}dx.
\]
Using Rodrigues' formula \cite[Eqn.\ (4.7.12)]{Szego-75-1} for $\ultra
{\lambda_n} \ell (x)$ in the equation above and employing the
duplication formula and other standard properties of the Gamma
function, one can obtain this expression:
\[
\hat\phi(\ell)=\frac{(-1)^\ell \omega_n
  \Gamma(\lambda_n+1)}{2^\ell\sqrt{\pi} \Gamma(\ell +\lambda_n
  +\frac12)}\int_{-1}^1\phi(x)\frac{d^\ell}{dx^\ell}
\left\{(1-x^2)^{\ell+\lambda_n-\frac12}\right\}dx,
\]
which holds for all $\ell$, even when $n=1$ -- i.e., $\lambda_1=0$.

Schoenberg \cite{Schoenberg-42-1} defined $\phi$ to be positive
definite if for every set of centers $X$ the matrix
$[\phi(\xi_j\cdot\xi_k)]$ is positive semidefinite. He showed that
$\phi$ is positive definite if and only if the coefficients satisfy
$\hat\phi(\ell)\ge 0$ for all $\ell$ and $\sum_{\ell=0}^\infty
\hat\phi(\ell)d_\ell <\infty$. If in addition $\hat\phi(\ell)>0$, then
$[\phi(\xi_j\cdot\xi_k)]$ is a positive definite matrix and one can
use shifts of $\phi$ to interpolate any function $f\in C(\sph^n)$ on
$X$. We will say that $\phi$ is a \emph{spherical basis function}
(SBF) in this case.

One usually makes the assumption that the sum $\sum_{\ell=0}^\infty
\hat\phi(\ell)d_\ell <\infty$, for then $\phi$ is continuous and
$\phi(1)=\|\phi\|_{L^\infty}$.  This is essential if we are doing
standard interpolation of a function from its values on $X$. However,
we are more interested in approximation than interpolation, and so we
will \emph{not} make this assumption here. Indeed, we will say that
any distribution $\phi$ for which $\hat\phi(\ell)>0$ for all $\ell$ is
a spherical basis function.  In general, we will be interested in SBFs
in $L^p$.

Zonal functions that satisfy $\hat\phi(\ell)>0$ for $\ell\ge L>0$ are
said to be \emph{conditionally positive definite} SBFs. In the RBF
theory on Euclidean space, the difference between strictly positive
definite RBFs and conditionally strictly positive definite RBFs is
significant. On $\sph^n$, this difference is less important: a
conditionally positive definite SBF differs from an SBF by a
polynomial of degree $L-1$. This does play a role in interpolation,
but is much less significant in approximation problems. That being the
case, unless there is a genuine need to distinguish between the two,
we will refer to both as simply SBFs.

Below we will list Fourier-Legendre expansion coefficients for some of the more
significant SBFs. Apart from certain Green's functions that we will do first, these SBFs are restrictions of Euclidean RBFs in $\RR^{n+1}$ to the $\sph^n$, which are themselves SBFs \cite[Corollary~4.3]{Narcowich-Ward-02-1}. These include Gaussians, multiquadrics, thin-plate splines, and Wendland functions. Such SBFs are RBFs expressed in terms of the Euclidean distance between $\xi$ and $\eta$ or its square, $\|\xi -\eta\|^2= 2 - 2\xi\cdot \eta$  and, with $t=\xi\cdot \eta$, these give rise to functions of $1-t$. 

\paragraph{Green's functions} Let $\beta>0$. The Green's function
solution to $\sfl_n^\beta G_\beta =\delta$ is a kernel with an expansion
in spherical harmonics having coefficients $\widehat
G_\beta(\ell,m)=(\ell+\lambda)^{-\beta}$. Properties of Green's
functions are discussed in more detail in Proposition~\ref{perturbation_approx_prop}. We simply remark that the kernel $G_\beta$ is an SBF that is in $L^1(\sph^n)$ for all $\beta>0$. For us, $G_\beta$ will play a significant role. The SBFs we consider will generally be of two types:  $\phi=G_\beta+G_\beta\ast \psi$, where $\psi$ is an $L^1$ zonal function, or $\phi$ will be in $C^\infty$. The first type includes the thin-plate splines and Wendland functions, and the second, the Gaussians and multiquadrics. 

\paragraph{Thin-plate splines} The thin-plate splines are defined in
\cite[Section 8.3]{Wendland-05-1}; their Fourier-Legendre coefficients
are found in \cite[\S 4.2]{Narcowich-etal-07-2}. These are given
below.
\begin{equation}\label{TPS}
\left.
\begin{array}{l}
\displaystyle{
\phi_s(t)=
\left\{
\begin{array}{cc}
  (-1)^{\lceil (s)_+\rceil}(1-t)^s, & s > - \frac{n}{2},  \ s \not\in \NN\\[5pt]
  (-1)^{s+1}(1-t)^s\log(1-t), &  s\in\NN.
\end{array}
\right.}\\[18pt]
  \hat\phi_s(\ell)=
  C_{s,n}\frac{\Gamma(\ell-s)}{\Gamma(\ell+s+n)}. 
\end{array}
\right\}
\end{equation}
where the factor $C_{s,n}$ is given by
\[
C_{s,n}:=2^{s+n}\pi^{\frac{n}{2}}\Gamma(s+1)\Gamma(s+\frac{n}{2}) 
\left\{
\begin{array}{ll}
 \frac{\sin(\pi s)}{\pi} & s > - \frac{n}{2},  \ s \not\in \NN\\[5pt]
1, & s\in\NN.
\end{array}
\right.
\]

Let $\nu=\ell+\lambda_n$. For large $\nu$, the Fourier-Legendre coefficients
$\phi_s(\ell)$ for the thin-plate splines have the asymptotic form
\begin{equation}
\label{tps_asymtotics}
\hat\phi_s(\ell)=C_{s,n}\nu^{-2s-n}\left(1+\sum_{j=1}^{p-1}
  G_{j}(n,s)\nu^{-j}+R_p(n,s,\nu)\right),
\end{equation}
where $R_p(n,s,\nu)=\calo(\nu^{-p})$ and $G_{j}(n,s)$ are defined in
\cite[p.~119]{Olver-74-1}. 

Two remarks. First, we have made use of $G_0(n,s)=1$ in the expansion
from \cite[p.~119]{Olver-74-1}. Second, when $s$ is an integer or
half-integer, $\hat\phi_s(\ell)$ is a rational function of $\ell$,
and, hence, of $\nu$. In that case, it follows that the series for
$\hat\phi_s(\ell)$ is actually a convergent power series in
$\nu^{-1}$. For other $s$, the expansion is only asymptotic.

From the structure of the expansions above and the properties of Green's functions listed in Proposition~\ref{perturbation_approx_prop}, we see that any finite linear combination of thin-plate splines
\begin{equation}
\label{lin_comb_tps}
\phi= \sum_{j=1}^m A_j \phi_{s_j}, \ - \frac{n}{2} < s_1<s_2 <\cdots <s_m,
\end{equation}
has the form
\begin{equation}
\label{lin_comb_tps_greens_fnc}
\phi=A_1(G_{2s+n}+G_{2s+n}\ast \psi), \ \psi\in L^1.
\end{equation}

\paragraph{Wendland functions} All of the SBFs we have discussed so
far are related to RBFs stemming from completely monotonic
functions. These RBFs have the property that they are strictly
positive definite or conditionally positive definite in $\RR^n$ for
all $n$. The corresponding SBFs are also positive definite in
$\sph^n$, again for all $n$. These RBFs are \emph{not} compactly
supported, however. This can be remedied, but there is a price: we
must give up positive definiteness beyond a certain dimension.

Wendland (cf. \cite[Section 9.4]{Wendland-05-1}) constructed families
of RBFs that are compactly supported on $0\le r\le R$, strictly positive
definite in Euclidean spaces of dimension $d$ or less, have smoothness
$C^{2k}$, and, within their supports, are polynomials of degree
$\lfloor \frac{d}{2} \rfloor + 3k+1$. The quantities $d$, $k$, and $R$ are parameters and may be adjusted as needed.

Restricting the Wendland functions to $\sph^n$ just requires setting $r=\sqrt{2(1 - t)}$ and $R=\sqrt{2(1 - t_0)}$, where $-1<t_0\le t \le1$. We will denote these functions by $\phi_{d,k}(t)$. The support of $\phi_{d,k}$ on  $\sph^n$ is then $0\le \theta\le \cos^{-1}(t_0)<\pi$. From \cite[Theorems 9.12 \& 9.13]{Wendland-05-1}, if $t>t_0$, then these functions are polynomials in $\sqrt{1-t}$ that may be put into the form,
\[
p_{d,k}(t) = e_1(1-t) + (1-t)^{k+\frac12}e_2(1-t),
\]
where $e_1$ and $e_2$ are polynomials having with $\deg e_1=\lfloor\frac12( \lfloor \frac{d}{2} \rfloor + 3k+1)\rfloor$ and $\deg e_2 = \lfloor \frac12(\lfloor \frac{d}{2} \rfloor + k)\rfloor$. Outside of this interval, the $\phi_{d,k}$ are identically $0$. Using a power series argument, we have that, near $t\gtrapprox t_0$, $\phi_{d,k}(t) = A(t-t_0)^{\lfloor \frac{d}{2} \rfloor + 2k+1}\big(1+ \calo(t-t_0)\big)$, from which it follows that $\phi_{d,k}(t) $ is piecewise $C^{\lfloor \frac{d}{2} \rfloor + 2k+1}$ near $t_0$. In addition, it follows that $\psi_{d,k}(t) := \phi_{d,k}(t) - p_{d,k}(t) $ is piecewise $C^{\lfloor \frac{d}{2} \rfloor + 2k+1}$ on the whole interval $[-1,1]$. Putting all of this together, we conclude that
\begin{equation}
\label{WendSBF-structure}
\phi_{d,k}(t) = e_1(1-t) + (1-t)^{k+\frac12}e_2(1-t) + \psi_{d,k}(t).
\end{equation}

Our aim is to use this decomposition to obtain large $\ell$ asymptotics for the Fourier-Legendre coefficients $\hat\phi_{d,k}(\ell)$ in $\sph^n$. This we now do.

\begin{prop}
\label{WendSBF-coef}
Let $m=\lfloor \frac{d}{2} \rfloor + 2k+1$. If $\ell>\deg e_1$, then
\[
\hat \phi_{d,k}(\ell) = (\ell+\lambda_n)^{-(2k+1+n)}\left(A_0+\frac{A_1}{\ell+\lambda_n}+ \calo(\ell+\lambda_n)^{-2} \right)+  \frac{\widehat{\sfl^m\psi_{d,k}}(\ell)}{(\ell+\lambda_n)^m}.
\]
Moreover, if we choose $\lfloor \frac{d}{2} \rfloor>n$, then the $\phi_{d,k}$ have the structure 
\[
\phi_{d,k} = \mbox{\rm polynomial } + A_0 \left(G_{2k+n}+G_{2k+n}\ast \tilde\psi\right), \ \tilde\psi \in L^.
\]
\end{prop}

\begin{proof}
The polynomial term $e_1(1-t)$  doesn't contribute to coefficients with $\ell>\deg e_1$. The term $(1-t)^{k+\frac12}e_2(1-t)$ is a linear combination of thin-plate splines, starting with $s=k+\frac12$. Thus it contributes the first term on the right above. By Remark~\ref{dom_L}, the function $\psi_{d,k}$ is in $H^2_m$, so it can be written as $\psi_{d,k}=\sfl_n^{-m}\sfl_n^m \psi_{d,k}$. The second term on the right follows directly from this fact. Finally, the form of the $\hat \phi_{d,k}(\ell)$'s  leads to the second statement.
\end{proof}

Before leaving the topic, we point out that, when $\lfloor \frac{d}{2} \rfloor>n$, we have determined the precise asymptotics of the Fourier-Legendre coefficients for the Wendland functions. Heretofore only upper and lower bounds were known.

\paragraph{Gaussians} The Fourier-Legendre coefficients for the
Gaussians, which are given below, may be found in \cite[Ex.\ 37, p.\
383]{Whittaker-Watson-65-1}, \cite[Example 5.2]{Mhaskar-etal-99-1}, and
\cite[\S 4.3]{Narcowich-etal-07-2}.
\begin{equation}
\label{gaussians}
\left.
\begin{array}{l}
  \gamma_\sigma(t)=e^{-2\sigma(1-t)}, \ \sigma>0, \\[6pt]
  \hat \gamma_\sigma(\ell) =
  2\pi\left(\frac{2\pi}{\sigma}\right)^{\lambda_n} e^{-\sigma} 
  I_{\lambda_n+\ell}(\sigma),
\end{array}
\right\}
\end{equation}
where $I_{\lambda_n+\ell}$ is an order $\lambda_n+\ell$ modified
Bessel function of the first kind. For all $\ell\ge 0$, the
coefficient $\hat \gamma_\sigma(\ell)$ satisfies this bound:
\cite[Proposition~4.3]{Narcowich-etal-07-2}:
\begin{equation}
\label{gaussian_bounds}
\frac{2 \sigma^\ell e^{-2\sigma}\pi^{\frac{n+1}2}
}{\Gamma(\ell+\frac{n+1}2)}\le \hat\gamma_\sigma(\ell) \le
\frac{2\sigma^\ell \pi^{\frac{n+1}2}}{\Gamma(\ell+\frac{n+1}2)}.
\end{equation}

\paragraph{Multiquadrics} The Hardy multiquadrics are treated in
\cite[\S 5]{Narcowich-etal-07-2}. The results are:
\begin{equation}
\label{multiquadric}
\left.
\begin{array}{l}
\mq_\alpha(t)=-\sqrt{\delta^2+2(1-t)},\ \delta>0. \\[12pt]
\begin{aligned}
  \widehat{\mq}_\delta(\ell) =&\frac{ \pi^{\lambda_n}
    \Gamma(\ell-1/2)}{(\alpha^2+2)^{\ell-1/2}
    \Gamma(\ell+\lambda_n+1)} \times \\
  &\hgeom{2}{1}\left(\frac{\ell-1/2}{2}, \frac{\ell+1/2}{2};
    \ell+\lambda_n+1; \frac{4}{(\delta^2+2)^2}\right).
\end{aligned}
\end{array}\right\}
\end{equation}
Here, $\hgeom{2}{1}$ is the usual hypergeometric function. Expressions
for Fourier-Legendre coefficients for generalized multiquadrics may be
found in \cite[\S 5]{Narcowich-etal-07-2}. Again, this time for $\ell$
sufficiently large, the coefficient $\widehat{\mq}_\delta(\ell)$
satisfies the following bound \cite[Proposition~5.1]{Narcowich-etal-07-2}:
\begin{equation}
\label{multiquadric_bounds}
C_1 \ell^{-\frac{n}{2}-1} \left( \frac{1}{\delta^2+2}
\right)^{\ell- \frac12} <\widehat{\mq}_\delta(\ell) <C_2
\ell^{-1-n} \left( \frac{2}{\delta^2+2} \right)^{\ell-\frac12},
\end{equation}

\paragraph{Ultraspherical generating functions} For $n\ge 2$, the
ultraspherical polynomials $\ultra {\lambda_n} \ell$ are frequently
defined in terms of the generating function \cite[Equation
(4.7.23)]{Szego-75-1} below:
\begin{equation}
\label{ultrasph_gen_fns}
\left.
\begin{array}{l}
u_{\lambda_n,w}(t) = (1 - 2tw+w^2)^{-\lambda_n},\ 1>w>0,\ n\ge 2 \\[6pt]
\hat u_{\lambda_n,w}(\ell) = w^\ell
\end{array}
\right\}
\end{equation}
When $n=1$, $\lambda_1=0$, the expansion is in terms the $T_\ell(t)$'s, the
Chebyshev polynomials of the first kind. In this case, the gerating
function is simply the Poisson kernel.
\begin{equation}
\label{poisson_kernel}
\left.
\begin{array}{l}
P_w(t) = \frac{1-w^2}{1 - 2tw+w^2},\ 1>w>0, \\[6pt]
\widehat P_w(\ell) =
\left\{ 
\begin{array}{l}
1, \ \ell=0,\\
2w^\ell, \ \ell\ge 1
\end{array}
\right.
\end{array}
\right\}
\end{equation}

\section{Approximation}\label{approximation}

The approximation part of the analysis makes use of kernels and
frames, which are related to them. These were studied in
\cite{Dai-07-1, Mhaskar-05-2, Mhaskar-etal-00-1,
  Mhaskar-Prestine-05-1,Narcowich-etal-06-1} and further developed in
\cite{Petrushev-Xu-2008-1}; we review them here, along with a number
of other results important to attaining the goals of this
paper. First, we will develop various types Marcinkiewicz-Zygmund
inequalities for the sphere. Although some of these were previously
derived
\cite{Mhaskar-etal-01-1,Mhaskar-etal-01-2,Narcowich-etal-06-1}, those
pertinent to both the approximation and stability analysis are new.

Second, using frames we establish a Bernstein inequality for spherical
polynomials. Again, using frames we establish various distance
estimates for $\phi\in H^1_\beta$ and we discuss Green's function
solutions to $\sfl_n^\beta G_\beta=\delta$. As we have mentioned
earlier, these form a very important class of SBFs.  Finally, at the
end of the section we will complete the approximation part of the
analysis.

\subsection{Kernels}

Let $\kappa(t)\in C^k(\RR)$, with $k\ge
\max\{2,n-1\}$, be even, not identically 0,  and satisfy
\begin{equation}
\label{kappa_condits}
|\kappa^{(r)}(t)| \le C_\kappa (1+|t|)^{r-\alpha}\ \mbox{for all }
 t\in \RR,\ r=0,\ldots,k,
\end{equation}
where $\alpha>n+k$ and $C_\kappa>0$ are fixed constants. We remark
that all compactly supported, $C^k$ functions that are even satisfy
(\ref{kappa_condits}). Functions in the Schwartz-class $\cals(\RR)$
that are even satisfy (\ref{kappa_condits}) for arbitrarily large $k$
and $\alpha$. Given such a $\kappa$, define the family of operators
\[
\sfk_{\varepsilon,n} := \kappa(\varepsilon \sfl_n) =
\sum_{\ell=0}^\infty \kappa(\varepsilon (\ell+\lambda_n))\sfp_\ell, \
0<\varepsilon\le 1,
\]
along with the associated family of kernels
\begin{equation}\label{K_kernel_def}
K_{\varepsilon,n}(\underbrace{\xi\cdot \eta}_{\cos\theta}) :=
\left\{
\begin{array}{ll}
\frac{1}{2\pi}\kappa(0) + 
\frac{1}{\pi}\sum_{\ell=1}^\infty
\kappa(\varepsilon\ell)\cos \ell \theta, &n=1,\\ [5pt]
\sum_{\ell=0}^\infty
\kappa(\varepsilon (\ell+\lambda_n))
\frac{\ell+\lambda_n}{\lambda_n\omega_n}\ultra
{\lambda_n} \ell (\cos \theta ), &n\ge 2,
\end{array}\right.
\end{equation}
where $\cos \theta =\xi\cdot \eta$ and $0<\varepsilon\le 1$.  

It is worthwhile noting that $\kappa(t) = e^{-t^2}$ satisfies
(\ref{kappa_condits}) and that the corresponding kernel is essentially
the heat kernel for $\sph^n$.

We will need several results concerning these kernels and
operators. First of all, we require the estimates on the $L^p$ norms
for the kernels. Material closely connected to the theorem below
appeared in \cite[Proposition 4.1]{Mhaskar-05-2}.

\begin{theorem}[{\cite[Theorem 3.5 \& Corollary 3.6]{Narcowich-etal-06-1}}]
\label{K_kernel_main_bnd}
Let $\kappa$ satisfy {\rm (\ref{kappa_condits})}, with $k\ge
\max\{2,n-1\}$.  If $0\le \theta\le \pi$, then there is a constant
$\beta_{n,k,\kappa}>0$ such the kernel $K_{\varepsilon,n}$ satisfies
the bound
\begin{equation}
\label{K_kernel_bnd_final}
|K_{\varepsilon,n}(\cos\theta)|\le \frac{\beta_{n,k,\kappa}}
{1+(\frac{\theta}{\varepsilon})^k} \varepsilon^{-n}.
\end{equation}
Moreover, we have that
\begin{equation}
\label{K_kernel_Lp_bound}
\|K_{\varepsilon,n}\|_p:=\|K_{\varepsilon,n}(\cos\theta)\|_{L^p(\sph^n)} \le 
C_{n,k,\kappa}\varepsilon^{-n/p'}.
\end{equation}
\end{theorem}

These operators can be applied to functions in $L^p(\sph^n)$ or even
distributions in $ \cald'(\sph^n)$, provided $\kappa$ decays fast
enough -- compact support will certainly work. As the result below
shows, all them are bounded operators taking $L^p(\sph^n)\to
L^q(\sph^n)$.

\begin{theorem}[{\cite[Theorem 3.7]{Narcowich-etal-06-1}}]
\label{K_kernel_op_bnd}
  If $\kappa$ satisfies (\ref{kappa_condits}), with $k>\max\{2,n\}$,
  then, for all $1\le p\le \infty$ and $1\le q\le \infty$, the
  operator $\sfk_{\varepsilon,n}\colon L^p(\sph^n)\to L^q(\sph^n)$ is
  bounded and its norm satisfies
\[
\|\sfk_{\varepsilon,n}\|_{p,q} \le
C_{n,k,\kappa}(4\omega_{n-1}
\varepsilon^n)^{-(\frac{1}{p}-\frac{1}{q})_+}\,,
\]
where $C_{n,k,\kappa}$ is a constant that depends only on $n, k, \kappa$, and
where $(x)_+=x$ for $x>0$ and $(x)_+= 0$ otherwise.
\end{theorem}

We point out that more can be said when $\kappa$ has restrictions on its support. The result below follows from the spherical harmonics of degree $L\sim1/\varepsilon$ being in the kernel of $\sfk_{\varepsilon,n}$ when $\kappa(t) = 0$ near $t=0$.
\begin{remark}
\label{kernel_support_condit}
If $\kappa(t)=0$ for $|t|\le 1$, then for any spherical harmonic in $\Pi_{L_\varepsilon}$, where $L_\varepsilon = \lfloor \varepsilon^{-1} - \lambda_n^{-1} \rfloor \sim \varepsilon^{-1}$ or less, then we have $g_\varepsilon:=\sfk_{\varepsilon,n}g = \sfk_{\varepsilon,n}(g-P)$, and hence $\| g_\varepsilon\|_q \le \|\sfk_{\varepsilon,n}\|_{p,q} E_{L_\varepsilon }(g)_p$.
\end{remark}

Another important result for $\kappa$ supported away from $t=0$ and having fast decay is the one below, which follows directly from Theorems~\ref{K_kernel_main_bnd}~and~\ref{K_kernel_op_bnd}. To simplify matters, we will assume that $\kappa$ is also compactly supported.

\begin{cor} \label{kernel_scaling_prop}
Let $k>\max\{2,n\}$. If the support of $\kappa$ is compact and does \emph{not} include $t=0$, then, for every fixed $\gamma$ in $\CC$, the function $\tilde\kappa(t):=|t|^\gamma \kappa(t)$ is also an even $C^k$ function that satisfies (\ref{K_kernel_def}).  Moreover, $\sfl^\gamma \sfk_{\varepsilon,n} =\varepsilon^{-\gamma} \tilde \sfk_{\varepsilon,n}.$ Finally, for real $\gamma$, we  have the two bounds below:
\begin{align*}
\|\sfl^\gamma \sfk_{\varepsilon,n}\|_{p,q} &\le C_{n,k,\tilde\kappa}
(4\omega_{n-1})^{-(\frac{1}{p}-\frac{1}{q})_+} \varepsilon^{-\gamma -n(\frac{1}{p}-\frac{1}{q})_+} \\
\|\sfl^\gamma \sfk_{\varepsilon,n}\delta\|_p &\le C_{n,k,\tilde\kappa}\varepsilon^{-\gamma-n/p'} ,
\end{align*}
where $\delta$ is the Dirac distribution and thus $\sfl^\gamma \sfk_{\varepsilon,n}\delta$ is the kernel for $\sfl^\gamma \sfk_{\varepsilon,n}$.
\end{cor}

\subsection{Marcinkiewicz-Zygmund inequalities}

Marcinkiewicz-Zygmund (MZ) inequalities provide equivalences between
norms defined through integrals and ones defined through discrete
sums. For $\sph^n$, these were developed in
\cite{Mhaskar-etal-01-1,Mhaskar-etal-01-2,Narcowich-etal-06-1}.  We
will need to adapt these MZ inequalities to estimate certain sums.

Let $X\subset \sph^n$ be the set of centers; also, let $q=q_X$,
$h=h_X$, and $\rho=\rho_X:=h/q$ be the separation radius, mesh norm,
and mesh ratio, respectively. We will need a decomposition of the
sphere into a finite number of non-overlapping regions. The Voronoi
tessellation corresponding to $X$ will serve our purpose here,
although many other decompositions will work as well.

Let $R_\xi$ be the Voronoi region containing $X$. Denote the
collection of these regions by $\calx =\{R_\xi\subset \sph^n\,|\,
\xi\in X\}$ and its partition norm by $\|\calx\|=\max_{\xi\in
  X}\{\mbox{diam}(R_\xi)\}$. It is easy to show that the following
geometric inequalities hold:
\begin{equation}
\label{q_h_calx}
h\le \|\calx\| \le 2h \ \mbox{and } \min_{\xi\in X}\mu(R_\xi) \ge c_n q^n.
\end{equation}
Here $c_n$ is a constant related to the volume of $\sph^n$. We will
need these later. For a sequence space version of results below, see
\cite[Proposition 4.1]{Mhaskar-06-1}.

\begin{prop}
\label{MZ_kernel}
Fix $\zeta\in\sph^n$ and $k\ge n+2$. Let
$K_\varepsilon(\eta):=K_{\varepsilon,n}(\eta\cdot\zeta)$. Then, there
is a constant $C=C_{n,\kappa,k}$ for which
\begin{equation}
\label{kernel_norm_est}
\bigg| \|K_{\varepsilon}\|_1 - \sum_{\xi\in
  X}\mu(R_\xi)|K_\varepsilon(\xi)| 
\bigg| \le C\left\{
\begin{array}{cc}
\| \calx\|/\varepsilon&  \|\calx\| \le \varepsilon  \\[5pt]
\left(\|\calx\|/\varepsilon\right)^n & \|\calx\| \ge  \varepsilon.
\end{array}
\right. 
\end{equation}
Moreover, if $\zeta\in X$, then 
\begin{equation}
\label{kernel_region_est}
\bigg| \int_{\sph^n-R_\zeta}|K_\varepsilon(\eta)|d\mu(\eta) - 
\sum_{X\ni\xi\ne \zeta}\mu(R_\xi)|K_\varepsilon(\xi)| \bigg| 
\le C_{n,\kappa,k} \left\{
\begin{array}{cc}
(\| \calx\|/\varepsilon)^{-1}&  \|\calx\| \le \varepsilon  \\[5pt]
\left(\varepsilon/\|\calx\| \right)^{k-n-2} & \|\calx\| \ge  \varepsilon.
\end{array}
\right. 
\end{equation}
\end{prop}

\begin{proof}
  The proof follows along the the lines of the one for
  \cite[Proposition~4.1]{Narcowich-etal-06-1}. Therefore, we will only
  sketch it here, referring the reader to \cite{Narcowich-etal-06-1}
  for the technical details.

  The inequalities in both (\ref{kernel_norm_est}) and
  (\ref{kernel_region_est}) involve bounding sums of contributions
  from each $R_\xi$ having the form
\[
D_\xi:=\bigg|\int_{R_\xi}|K_\varepsilon(\eta)|d\mu(\eta) -
\mu(R_\xi)|K_\varepsilon(\xi)| \bigg| \le
\int_{R_\xi}|K_\varepsilon(\eta)-K_\varepsilon(\xi)|d\mu(\eta).
\]
Take $\zeta$ to be the north pole of the sphere and $\theta$ to be the
co-latitude. Divide the sphere into $M \sim \pi/\| \calx\| $ bands,
$B_m$, in which $(m-1)\pi/M \le \theta \le m\pi/M$, $m=1,\ldots, M$.
Each $R_\xi$ can have non-trivial intersection with at most two
adjacent bands, because $\mbox{diam}(R_\xi) \le \|\calx\| \sim
\pi/M$. Thus, if $R_\xi\subset B_m\cup B_{m+1}$, then its lowest and
highest co-latitudes satisfy $(m-1)\pi/M\le \theta_\xi^- \le
\theta_\xi^+\le (m+1)\pi/M$. As is shown in
\cite{Narcowich-etal-06-1}, for $m=2,\ldots,M-1$ the sum of the
$D_\xi$ from all $R_\xi \subset B_m\cup B_{m+1}$ is bounded above by
the quantity,
\begin{equation}
\label{D_non_cap_bnd}
\sum_{R_\xi \subset B_m\cup B_{m+1}} D_\xi  \le  
\frac{C_{n,\kappa,k}}{M\varepsilon}
\int_{\frac{m-1}{M\varepsilon}\pi}^{\frac{m+1}{M\varepsilon}\pi}
\frac{t^n}{1+t^k} dt .
\end{equation}
If $R_\xi \ni \zeta$, then dealing with the corresponding $D_\xi$ can
be done by estimating the integral that bounds the contribution from
the region $R_\xi$ in the cap $0\le \theta\le 2\pi/M$,
\begin{equation}
\label{D_cap_bnd}
D_\xi \le C'_{n,\kappa,k}(M\varepsilon)^{-n} 
\int_0^{\frac{2\pi}{M\varepsilon} }\frac{ t dt  }
{1+t^k}\le \frac{C''_{n,\kappa,k}}{(M\varepsilon)^n} \left\{
\begin{array}{cc}
(M\varepsilon)^{-2} & M\varepsilon\ge 1 \\ [5pt]
1 & M\varepsilon \le 1 .
\end{array}
\right. 
\end{equation}
Now, let $M=\lfloor \pi/\|\calx\|\rfloor$, precisely. Adding up the
$D_\xi$ for all $\xi\in X$ yields the bound in
(\ref{kernel_norm_est}), which was implicit in the proof of
\cite[Proposition~4.1]{Narcowich-etal-06-1}.

To get (\ref{kernel_region_est}), we need to adjust $M$ so that all
$R_\xi \not\ni \zeta$ are contained in the bands $B_m\cup B_{m+1}$,
$m=2,\ldots,M-1$. This is easy to do. Just take $M=\lfloor
(\pi-q)/\|\calx\|\rfloor$. Summing the $D_\xi$ bounded in
(\ref{D_non_cap_bnd}) and taking care of some double counting yields
\begin{eqnarray*}
  \bigg| \int_{\sph^n-R_\zeta}|K_\varepsilon(\eta)|d\mu(\eta) -
  \sum_{X\ni\xi\ne \zeta}\mu(R_\xi)|K_\varepsilon(\xi)| \bigg| &\le& 
  \frac{C_{n,\kappa,k}}{M\varepsilon}
  \int_{\frac{\pi}{M\varepsilon}}^{\frac{\pi}{\varepsilon}}
  \frac{t^n}{1+t^k} dt \\
  &\le& \frac{C_{n,\kappa,k}}{M\varepsilon}
  \int_{\frac{\pi}{M\varepsilon}}^{\infty}
  \frac{t^n}{1+t^k} dt \\
  &\le & C_{n,\kappa,k}\left\{
\begin{array}{cc}
(M\varepsilon)^{-1} & M\varepsilon\ge 1 \\ [5pt]
(M\varepsilon)^{k-n-2} & M\varepsilon \le 1 ,
\end{array}
\right. 
\end{eqnarray*}
from which (\ref{kernel_region_est}) follows easily.
\end{proof}

Let $f\in L^1(\sph^n)$ and set $f_\varepsilon := K_{\varepsilon,n}*f$;
the function $f$ is \emph{not} assumed to be zonal. We wish to
estimate the difference $ E_\calx:=\big|\|f_\varepsilon\|_1 -
\sum_{\xi\in X}|f_\varepsilon(\xi)|\mu(R_\xi)\big| $. It is
straightforward to show that
\[
E_\calx \le \sum_{\xi\in X}\int_{R_\xi}|f_\varepsilon(\eta) -
f_\varepsilon(\xi)|d\mu(\eta) \le \sup_{\zeta\in\sph^n}
F_{\varepsilon,\calx}(\zeta) \|f\|_1\,,
\]
where $F_{\varepsilon,\,\calx}(\zeta):= \sum_{\xi\in
  X}\int_{R_\xi}\big| K_{\varepsilon,n}(\eta\cdot\zeta) -
K_{\varepsilon,n} (\xi\cdot\zeta)\big|d\mu(\eta)$, which is the
quantity estimated in Proposition~\ref{MZ_kernel}. Applying
that proposition and Remark~\ref{kernel_support_condit}, we
obtain the desired estimate below.

\begin{cor}\label{MZ_function_cor}
  Let $\kappa$ satisfy (\ref{kappa_condits}), with $k\ge n+2$, and,
  for $f\in L^1(\sph^n)$, let $f_\varepsilon=K_{\varepsilon,n}*f$. If
  $\calx$ is the decomposition of $\sph^n$ described above,
  $\|\calx\|\ge \varepsilon$ and $L_\varepsilon = \lfloor
  \varepsilon^{-1} - \lambda_n^{-1} \rfloor \sim
  \varepsilon^{-1}$. then
\begin{equation}
\label{MZ_function_est}
\bigg|\|f_\varepsilon\|_1 - \sum_{\xi\in
  X}|f_\varepsilon(\xi)|\mu(R_\xi)\bigg| 
\le C_{n,\kappa,k} \left(\|\calx\|/\varepsilon\right)^n \left\{
\begin{array}{cc}
E_{L_\varepsilon}(f)_1,& \kappa(t) = 0, \, |t|\le 1.  \\[4pt]
\|f\|_1, & \mbox{otherwise},
\end{array}
\right.
\end{equation}
\end{cor}

\begin{remark}\label{MZ_zonal_fnct_est}
  If $f$ is zonal, i.e. $f(\xi) = \psi(\xi\cdot \zeta)$, then the
  right side (\ref{MZ_function_est}) is independent of the variable
  $\zeta$. Also, the strict inequality $\|\calx\|\ge \varepsilon$ isn't absolutely necessary. The results still hold when $\|\calx\|$ and $\varepsilon$ are comparable.
\end{remark}

For the most part, we will use these results to bound the sums
$\big|\sum_{\xi\in X}a_\xi f_\varepsilon(\xi)\big|$, under the
assumption that $ \|\calx\| \ge \varepsilon$. Using
Corollary~\ref{MZ_function_cor} for that case, we see that
\begin{eqnarray}
  \bigg|\sum_{\xi\in X}a_\xi f_\varepsilon(\xi)\bigg| &\le& 
  \frac{|a|_\infty}{\min_{\xi\in X}\mu(R_\xi)}\sum_{\xi\in X} 
  \mu(R_\xi)|f_\varepsilon(\xi)|  \nonumber \\
  & \le&  \frac{|a|_\infty}{\min_{\xi\in X}\mu(R_\xi)}\left(\|f_\varepsilon \|_{L^1}
    + C_{n,\kappa,k}\left(\|\calx\|/\varepsilon\right)^n \|f\|_{L^1}
  \right)
  \nonumber
\end{eqnarray}
From Theorem~\ref{K_kernel_op_bnd}, (\ref{q_h_calx}), and $h=\rho q$,
with $L_\varepsilon \sim \varepsilon^{-1}$ and $\rho q \approx
\|\calx\|\ge \varepsilon$. we have that
\begin{equation}
\label{SUM_function_est}
\bigg|\sum_{\xi\in X}a_\xi f_\varepsilon(\xi)\bigg|  \le C \rho^n
\varepsilon^{-n} |a|_\infty  \left\{
\begin{array}{cc}
E_{L_\varepsilon}(f)_1,& \mbox{if } \kappa(t) = 0, \, |t|\le 1,  \\[4pt]
\|f\|_1, & \mbox{otherwise},
\end{array}
\right.
\end{equation}
If $f$ is a zonal function, then, by Remark~\ref{MZ_zonal_fnct_est}, we may use the $\|\cdot\|_\infty$ norm on the left above.

We want to make the same kind of estimate, but for $f$ being replaced
by $\delta_\zeta$, the usual Dirac delta function. Thus
$f_\varepsilon$ is replaced by
$K_\varepsilon(\cdot):=K_{\varepsilon,n}*\delta(\cdot)=K_{\varepsilon,n}((\cdot)\cdot\zeta)$. A
nearly identical argument to the one used above, coupled with
(\ref{kernel_norm_est}) for $\|\calx\|\ge \varepsilon$ and the bound
on $\|K_\varepsilon\|_1$ from Theorem~\ref{K_kernel_op_bnd}, results
in
\begin{equation}
\label{LC_kernel_bnd_pointwise}
\bigg|\sum_{\xi\in X}a_\xi K_{\varepsilon}((\cdot)\cdot \xi)\bigg|  \le C\rho^n\varepsilon^{-n} |a|_\infty .
\end{equation}
The constants on the right above hold uniformly, so we thus have
\begin{equation}
\label{LC_kernel_bnd}
\bigg\|\sum_{\xi\in X}a_\xi K_{\varepsilon}((\cdot)\cdot \xi)\bigg\|_\infty   \le C\rho^n\varepsilon^{-n} |a|_\infty. 
\end{equation}

The two bounds above are very similar and can be used in
combination. They will be needed to complete the approximation part of
the analysis. There is another bound, somewhat different from these
two, that we will need in section~\ref{stability}:

\begin{lemma} If $\rho q \sim \| \calx \| \ge \varepsilon>0$ and if $k\ge n+2$, then
\begin{equation}
\label{SUM_kernel_est}
\max_{\zeta\in X}\sum_{X\ni\xi\ne \zeta} |K_{\varepsilon,n}(\xi\cdot
\zeta)| \le C_{n,\kappa,k} q^{-n}.
\end{equation}
\end{lemma}

\begin{proof}
In equation (\ref{kernel_region_est}), Proposition~\ref{MZ_kernel}, again
for $\|\calx\|\ge \varepsilon$, an argument similar to the ones used above 
gives us
\begin{equation*}
\sum_{X\ni\xi\ne \zeta} |K_{\varepsilon,n}(\xi\cdot \zeta)| \le
C'_{n,\kappa,k} q^{-n} \int_{\sph^n-R_\zeta}|K_{\varepsilon,n}(\eta\cdot\zeta)|d\mu(\eta) +
C''_{n,\kappa,k}q^{-n}\left(\frac{\varepsilon}{\|\calx\|} \right)^{k-n-2} 
\end{equation*}
Using $\int_{\sph^n-R_\zeta}|K_{\varepsilon,n}(\eta\cdot\zeta)|d\mu(\eta)\le
\|K_{\varepsilon,n}\|_1\le C_{n,\kappa,k}$, $\frac{\varepsilon}{\|\calx\|}\le 1$, and maximizing over $\zeta\in X$, we obtain (\ref{SUM_kernel_est}).
\end{proof}

This estimate is more delicate than (\ref{LC_kernel_bnd}), because the term missing from
the sum is $K_{\varepsilon,n}(\zeta \cdot \zeta) =  K_{\varepsilon,n}(1)$, which turns out to be $\calo(\varepsilon^{-n})$. For $\varepsilon/q$ small enough, the sum  (\ref{SUM_kernel_est}) will be majorized by  $ K_{\varepsilon,n}(1)$. This is needed as part of a diagonal dominance argument.

\subsection{Frames}\label{frames}

We now address the question of the frame decomposition mentioned
previously. Our approach follows the one in
\cite{Narcowich-etal-06-1}. As mentioned earlier, others are certainly
possible. For this, we need a function $a \in C^k(\RR)$,which we may
assume is even, with support in $[-2,-\frac12]\cup [\frac12,2]$, and
satisfying $|a(t)|^2 + |a(2t)|^2 \equiv 1$ on $[\frac12,1]$. Such a
function can be easily constructed out of an \emph{orthogonal wavelet
  mask} $m_0$ \cite[\S 8.3]{Daubechies-92-1}. In fact, if $m_0(\xi)\in
C^{k+1}$, then $a(t) := m_0(\pi \log_2(|t|))$ on $[-2,-\frac12]\cup
[\frac12,2]$, and $0$ otherwise, is a $C^k$ function that satisfies
the appropriate criteria. Define $b\in C^k(\RR)$ by
\begin{equation}
\label{b_def}
b(t):=
\left\{
\begin{array}{cc}
1 & |t|\le 1 \\[3pt]
|a(t)|^2 & |t|>1.
\end{array}\right.
\end{equation}
Using the properties of $a$ we see that $ \sum_{j=-\infty}^J
|a(t/2^j)|^2 = b(t/2^J)$ if $t>0$. In the sum on the left, only terms
with $j\ge \lfloor \log_2(t) \rfloor$ contribute. Terms with
$j<\lfloor \log_2(t) \rfloor$ are identically 0.

The quantity $\lfloor \log_2(t) \rfloor$ is obviously important. On
the $\sph^n$, the integer that corresponds to it is this:
\begin{equation}\label{j_n_def}
j_n:=
\left\{
\begin{array}{cc}
0&n=1,\\
\lfloor \log_2(\lambda_n) \rfloor &n\ge 2.
\end{array}\right.
\end{equation}

The integer $j_n$ helps us in defining our frame operators, which we
now do. Let $\sfa_j := a(2^{-j-j_n}\sfl_n)$ and $\sfb_j :=
b(2^{-j-j_n}\sfl_n)$.  Taking into account the support of $a$, we have
$\sfb_J = \sum_{j=0}^J \sfa_j \sfa_j^\ast$ for $n\ge 2$ . For $n=1$, a
projection $\sfp_0$ onto the constant function enters, and $\sfb_J =
\sfp_0+\sum_{j=0}^J \sfa_j \sfa_j^\ast$. We will need the following
approximation result concerning these operators.

\begin{prop}[{\cite[Proposition
    5.1]{Narcowich-etal-06-1}}]\label{B_J_approx_thm}
  Let $k>\max\{n,2\}$, and let $b$ be defined by (\ref{b_def}), with
  $a\in C^k(\RR)$. If $f\in L^p(\sph^n)$, $1\le p\le \infty$, and if
  $L>0$ is an integer such that $2^{-J-j_n}\le(L+\lambda_n)^{-1}$,
  then
\begin{equation}
\label{B_J_approx_id}
\|f-\sfb_J f\|_p\le C_{b,k,n}
E_L(f)_p\,,\ E_L(f)_p:=\dist_{L^p}(f,\Pi_L).
\end{equation}
Also, for $1\le p<\infty$ or, if $p=\infty$, for $f\in
C(\sph^n)$, we have $\label{B_J__approx_id_lim}
\lim_{J\to\infty} \sfb_J f = f$.
\end{prop}

\paragraph{Bernstein/Nikolskii inequalities.} There are several
inequalities that follow easily using frames. We will give a
Nikolskii-type inequality, which is a well-known inequality (\cite
[Proposition 2.1]{Mhaskar-etal-99-1} and \cite[\S
3.5]{Narcowich-etal-06-1}), From our point of view, the most important
inequality derived here is a Bernstein theorem for spherical
polynomials \cite[Theorem 2 (Eng. transl.)]{Rustamov-92-1}. An
independent proof is given in \cite[Proposition
4.3]{LeGia-Mhaskar-06-1}. For the convenience of the reader, short
proofs for both are given below.
\begin{theorem}
\label{bernstein-nikolskii-poly-ineqs}
Let $S\in \Pi_L$. Then, for $1\le p,q \le \infty$ and for $\gamma>0$,
we have
\begin{eqnarray}
\label{nikolskii_poly_est}
\mbox{\rm (Nikolskii) }\ \| S\|_q  
&\le& C_{p,q,n}L^{n(\frac{1}{p}-\frac{1}{q})_+}\|S\|_p\\[5pt]
\label{bernstein_poly_est}
\mbox{\rm (Bernstein) }\  \|S\|_{H^p_\gamma} &\le& C_{n,\gamma} L^\gamma \|S\|_p
\end{eqnarray}
\end{theorem}
\begin{proof}
  Let $\gamma>0$ and suppose $L+\lambda_n\le 2^{J+j_n}$. From the
  definition of $\sfb_J$, it is easy to see that $\sfb_J$ reproduces
  $\Pi_L$, and so $\sfb_JS=S$ for all $S\in \Pi_L$.  By Theorem
  \ref{K_kernel_op_bnd}, with $\kappa=b$ and
  $\varepsilon=2^{-J-j_n}\sim L^{-1}$, we see that $\| S\|_q \le
  C_{p,q,n}L^{n(\frac{1}{p}-\frac{1}{q})_+}\|S\|_p, \ S\in
  \Pi_L$. Dependence of the constants on $b$ and $k$ disappears upon
  taking the infimum over these two quantities, yielding
  (\ref{nikolskii_poly_est}).

  We now establish the Bernstein inequality. If $S\in \Pi_L$, then so
  is $\sfl^\gamma S$, and we have that $\sfb_J \sfl_n^\gamma
  S=\sfl_n^\gamma S$, provided $L+\lambda_n\le 2^{J+j_n}$. Using the expansion $\sfb_J = \sum_{j=0}^J \sfa_j \sfa_j^\ast$, we see that
\[
\sfl^\gamma S = \sum_{j=0}^J \sfa_j \sfa_j^\ast \sfl^\gamma S =
\sum_{j=0}^J \sfl^\gamma\sfa_j \sfa_j^\ast S.
\]
Consequently, we have that $\|S\|_{H^\gamma_p} = \|\sfl^\gamma S\|_p
\le \sum_{j=0}^J \|\sfl^\gamma\sfa_j \sfa_j^\ast \|_{p,p}
\|S\|_p$. Applying Corollary~\ref{kernel_scaling_prop}, with
$\kappa(t)=|a(t)|^2 $ and $\varepsilon = 2^{-j-j_n}$ for each $j$,
then yields this:
\begin{eqnarray*}
  \|S\|_{H^\gamma_p} =&\le& \bigg( \sum_{j=0}^J
  2^{(j+j_n)\gamma}\bigg) 
  C_{a,n,\gamma} \| S\|_p \\
  &\le& \frac{2^{(J+j_n+1)\gamma}-2^{j_n\gamma}}{2^\gamma-1}
  C_{a,n,\gamma} 
\| S\|_p \le L^\gamma C_{a,n,\gamma}\|S\|_p\,,
\end{eqnarray*}
where again $L\sim 2^{J+j_n}$. In the last inequality of the chain
above, we can take the infimum over all $a$ satisfying the requisite
conditions. This yields (\ref{bernstein_poly_est})
\end{proof}

\paragraph{Distance estimates.}
Frames can be used to estimate the distance in $L^p(\sph^n)$
from the polynomials to a function in a smoother space. If $f\in L^p$,
let $E_L(f)_p:=\dist_{L^p}(f,\Pi_L)$. Because $\sfb_Jf$ is a spherical
polynomial in $\Pi_{2^{J+j_n+1}}$, we have
\[
E_L(f)_p \le \|f-\sfb_Jf\|_p, \ L+\lambda_n\le 2^{J+j_n+1}.
\]
And because $\sfb_Jf$ converges to $f$ in all $L^p$, $1\le p <\infty$
and $p=\infty$ if $f\in C(\sph^n)$, we also have that
\[
E_L(f)_p \le \|f-\sfb_Jf\|_p \le \sum_{j=J+1}^\infty \|\sfa_j\sfa_j^* f \|_p,
\]
where the right side above may be infinite. Now, suppose that
$f=\sfl_n^\gamma h$, $h \in H^q_\beta(\sph^n)$,  In that case, we
have $\sfa_j\sfa_j^* \sfl_n^\gamma h =
\sfl_n^{-(\beta-\gamma)}\sfa_j\sfa_j^* \sfl_n^\beta h$. From this
and Corollary~\ref{kernel_scaling_prop}, with $p\leftrightarrow q$, we arrive at
\[
\|\sfa_j\sfa_j^* \sfl_n^\gamma h \|_p = \| \sfl_n^{-(\beta-\gamma)}\sfa_j\sfa_j^* \sfl_n^\beta h\|_p \le 2^{-(\beta-\gamma -n(\frac{1}{q}-\frac{1}{p})_+)(j+j_n)}C_{n,k,a} \|h\|_{H^q_\beta}
\]
Insert this in the equation above, sum the appropriate geometric series, and take $L\sim
2^{J+j_n}$ to get
\[
E_{2^{J+j_n}}(\sfl_n^\gamma h)_p \le C'_{\beta-\gamma,a,k,n}2^{-(\beta-\gamma -n(\frac{1}{q}-\frac{1}{p})_+)(J+j_n)}\|h\|_{H^q_\beta}\,,
\]
which was essentially obtained by Kamzolov \cite{Kamzolov-82-1}. Now,
since the left side above is unchanged if we replace $\sfl_n^\gamma$ by $\sfl_n^\gamma -S$, $S\in \Pi_{2^{J+j_n}}$, we can
replace $\|h\|_{H^q_\beta}$ by $E_{2^{J+j_n}}(\sfl_n^\beta
h)_q$. Collecting these results yields the proposition below.

\begin{prop}\label{H^p_gamma_H^q_beta}
  Let $\gamma\ge 0$, and $\beta>\gamma+n(\frac{1}{q}-\frac{1}{p})_+)$, where
  $1\le p,q \le\infty$. If $h\in H^q_\beta$, then there is a constant $C=C_{n,\beta, \gamma,a}$ such that
\[
E_{2^{J+j_n}}(\sfl_n^\gamma h)_p \le
\|(I-\sfb_J)h\|_{H^p_\gamma} \le C_{n,\beta,\gamma,
  a}2^{-(\beta-\gamma-n(\frac{1}{q}-\frac{1}{p})_+))(J+j_n)} E_{2^{J+j_n}}(\sfl_n^\beta h)_q
\,.
\]
\end{prop}

\paragraph{Green's functions and their properties.} Let
$\beta>n/p'$. Recall that the Green's function solution to $\sfl_n^\beta
G_\beta =\delta$ is a kernel with an expansion in spherical harmonics
having coefficients $\widehat G_\beta(\ell,m)=(\ell+\lambda)^{-\beta}$. Properties of Green's
functions (pseudo-differential operator kernels, really) on manifolds
have been studied extensively (cf.\ \cite{Hormander-87-1}). Our aim here is to use frames to obtain properties and various distance estimates that we need here quickly, and in a self contained way, for SBFs of the form $\phi_\beta = G_\beta+G_\beta\ast \psi$, where $\psi\in L^1$. Because the $\phi_\beta$'s are not in any of the Bessel-Sobolev spaces $H^p_\beta$, they have to be treated separately from the class in Proposition~\ref{H^p_gamma_H^q_beta} above 

We will begin with Green's functions themselves. Note that $\sfa_j
\sfa_j^*G_\beta= \sfl_n^{-\beta} \sfa_j \sfa_j^*\delta$. Since $ \sfa_j
\sfa_j^*=|a|^2(2^{-j-j_n}\sfl_n)$, where both $a$ and, of course,
$|a|^2$, have compact support that excludes $t=0$. we may apply Corollary~\ref{kernel_scaling_prop}, with $\varepsilon_j:=2^{-(j+j_n)}$.

\begin{equation}
\label{p-j-est-green_fnct}
\|\sfa_j \sfa_j^*G_\beta \|_p \le  C_{n,\beta,a}
\varepsilon_j^{\beta - n/p'}= C_{n,\beta,a} 2^{-(\beta-n/p') (j+j_n)}.
\end{equation}
Thus, for $\beta>n/p'$, the terms in $\sum_{j=0}^\infty \sfa_j
\sfa_j^*G_\beta $ are bounded by a geometric series, and so the
Weierstrass $M$ test implies that the series converges in $L^p$. That
is, we have shown that when $\beta>n/p'$ the limit $ \lim_{J\to\infty}
\sfb_JG_\beta$ is in $L^p$. A simple duality argument then shows that
the kernel $G_\beta = \lim_{J\to\infty} \sfb_JG_\beta$ in
$L^p(\sph^n)$.  Summing the geometric series in
(\ref{p-j-est-green_fnct}) yields $\|G_\beta -\sfb_J G_\beta\|_p \le C\,
2^{-(\beta-n/p') (J+j_n)}$.

These results also give us error bounds in $H^p_\gamma(\sph^n)$. If
$\gamma\ge 0$, then $\sfl^\gamma G_\beta = G_{\beta-\gamma}$ and
$\sfl^\gamma \sfb_J G_\beta = \sfb_J G_{\beta - \gamma}$.  This and
the estimate above imply that if in addition $\beta>\gamma+n/p'$, then
\begin{equation}
\label{G-beta-H-gamma-dist-est}
\|G_\beta -\sfb_J G_\beta\|_{H^p_\gamma} = \|G_{\beta-\gamma} -\sfb_J 
G_{\beta-\gamma}\|_p \le  C\, 2^{-(\beta-\gamma-n/p')(J+j_n)}.
\end{equation}

Perturbations of $G_\beta$ can be dealt with, too. Let $\psi$ be in
$L^1$. By Theorem~\ref{K_kernel_op_bnd}, (\ref{p-j-est-green_fnct})
and Remark~\ref{kernel_support_condit}, we have that, for all $j\ge
J$,  
\[
\|\sfa_j \sfa_j^*G_\beta\ast \psi \|_p \le \|\sfa_j
\sfa_j^*G_\beta\|_{1, p\,}E_{2^{j+j_n}}(\psi)_1 \le C 2^{-(\beta-\gamma)(j+j_n)}
E_{2^{J+j_n}}(\psi)_1
\]
Summing a geometric series and using (\ref{G-beta-H-gamma-dist-est}),
we arrive at the following bound.

\begin{prop}\label{perturbation_approx_prop}
  Let $\gamma\ge 0$, $\beta>\gamma+n/p'$, $\varepsilon_j = 2^{-(j+j_n)}$, and
  let $\psi\in L^1$ be a zonal function. If
  $\phi_\beta=G_\beta + G_\beta\ast \psi$, then $\phi_\beta \in H^p_\gamma$
  and there is a constant $C=C_{n,\beta,\gamma,a}$, which depends only on $n$, $\beta$,
  $\gamma$, and the function $a$, such that
\begin{equation}
\label{perturbed_approx_est}
E_{2^{J+j_n}}(\sfl^\gamma\phi_\beta)_p \le \|(I - \sfb_J )\phi_\beta
\|_{H^p_\gamma} \le 
C_{n,\beta,\gamma,a}\bigg(1+\varepsilon_J^{n/p'}E_{2^{J+j_n}}(\psi)_1 \bigg)
\varepsilon_J^{\beta-\gamma-n/p'}.
\end{equation}
\end{prop}

\subsection{Approximation analysis}
\label{approx_analysis}

The task at hand is to estimate the norms
$\|(I-\sfb_J)g\|_{H^p_\gamma}/|a|_p$, where $g\in \calg_{X,\phi}$. Our
approach will be to carry this out for $p=1$ and $p=\infty$, then use
the Riesz-Thorin theorem to obtain the result for all intermediate
values of $p$.

The easier of the two cases is $p=1$. Since $g\in \calg_{X,\phi}$, we
have $g=\sum_{\xi\in X}a_\xi \phi((\cdot)\cdot \xi)$. Again, let $\varepsilon_j = 2^{-(j+j_n)}$. From the
triangle inequality, the rotational invariance of the norms involved,
and Proposition~\ref{H^p_gamma_H^q_beta} and
Proposition~\ref{perturbation_approx_prop} it follows that
\begin{align*}
\|(I-\sfb_J)g\|_{H^1_\gamma}&\le |a|_1 \| (I-\sfb_J)\phi \|_{
  H^1_\gamma} \\
&\le  C \varepsilon_J^{\beta-\gamma}|a|_1
\left\{
\begin{array}{cl}
E_{2^{J+j_n}}(\sfl_n^\beta \phi)_1 & \phi\in H^1_\beta\,,\\
(1+E_{2^{J+j_n}}(\psi)_1 ) & \phi=G_\beta + G_\beta\ast \psi\,.
\end{array}\right.
\end{align*}

The $p=\infty$ case requires using frames.  Again, we have that
\[
\|(I-\sfb_J)g\|_{H^\infty_\gamma}\le \sum_{j=J+1}^\infty
\|\sfa_j\sfa_j^*\sfl_n^\gamma g\|_\infty \,,
\]
where $ \sfa_j\sfa_j^*\sfl_n^\gamma g = \sum_{\xi\in X}a_\xi
\sfa_j\sfa_j^*\sfl_n^\gamma \phi((\cdot)\cdot \xi)$. By equation
(\ref{SUM_function_est}), with $f=\sfl_n^\gamma \phi$,
$K_{\varepsilon_j,n}$ corresponding to $\kappa(t)=|a(t)|^2$,
$h \ge \varepsilon_J\ge \varepsilon_j$, all $j\ge J$, and $L_\varepsilon\sim 2^{j+j_n}$, we have
\[
\|\sfa_j\sfa_j^*\sfl_n^\gamma g\|_\infty = \bigg\|\sum_{\xi\in X}a_\xi
f_\varepsilon((\cdot)\cdot \xi) \bigg\|_\infty \le C \rho^n
\varepsilon_j^{-n}|a|_\infty \,E_{2^{j+j_n}}(\sfl_n^\gamma \phi)_1\,.
\]
By Proposition~\ref{H^p_gamma_H^q_beta} and
Proposition~\ref{perturbation_approx_prop}, with $J$ there replaced by
$j$, $p=\infty$, we have
\[
\|\sfa_j\sfa_j^*\sfl_n^\gamma g\|_\infty \le C |a|_\infty \rho^n
\varepsilon_j^{\beta-\gamma-n}
\left\{
\begin{array}{cl}
E_{2^{j+j_n}}(\sfl_n^\beta \phi)_1 & \phi\in H^1_\beta\,,\\
(1+\varepsilon_j^n E_{2^{j+j_n}}(\psi)_1 ) & \phi=G_\beta + G_\beta\ast \psi\,.
\end{array}\right.
\]
Since $E_{2^{j+j_n}}(f)_1\le E_{2^{J+j_n}}(f)_1$ when $j\ge J$, in the
inequality above we may replace the distances with respect to
$2^{j+j_n}$ with ones with respect to $2^{J+j_n}$. Doing so and again
summing a geometric series, we obtain
\[
\|(I-\sfb_J)g\|_{H^\infty_\gamma}\le C |a|_\infty \rho^n
\varepsilon_J^{\beta-\gamma-n}
\left\{
\begin{array}{cl}
E_{2^{J+j_n}}(\sfl_n^\beta \phi)_1 & \phi\in H^1_\beta\,,\\
(1+\varepsilon_J^n E_{2^{J+j_n}}(\psi)_1 ) & \phi=G_\beta + G_\beta\ast \psi\,.
\end{array}\right.
\]
Applying the Riesz-Thorin theorem in conjunction with the bounds
above, we complete the approximation part of the problem:

\begin{theorem} 
\label{main_approximation_thm}
Let $\gamma\ge 0$, $1\le p\le \infty$, $\beta>\gamma+n/p'$, $\varepsilon_j = 2^{-(J+j_n)}$.  If $h_X\ge \varepsilon_j$ and if $g\in \calg_{X,\phi}$, then
\begin{equation}
\label{final_approx_est}
\frac{\|(I-\sfb_J)g\|_{H^p_\gamma}}{|a|_p} \le C \rho^{n/p'}
\varepsilon_J^{\beta-\gamma-n/p'}
\left\{
\begin{array}{cl}
E_{2^{J+j_n}}(\sfl_n^\beta \phi)_1 & \phi\in H^1_\beta\,,\\
(1+E_{2^{J+j_n}}(\psi)_1 ) & \phi=G_\beta + G_\beta\ast \psi\,.
\end{array}\right.
\end{equation}
\end{theorem}

\section{Stability}\label{stability}

The problem that we address here is estimating the norm $|a|_p$ in
terms of the $L^p(\sph^n)$ norm of $g$, where $g(\bfx)=\sum_{\xi\in
  X}a_\xi \phi(\bfx \cdot \xi)$ and $\phi\in L^p$ is an SBF.
Specifically, we wish to estimate the $p$-norm stability ratio
\[
\maxg_{\calg,\,p}:= \max_{\calg \ni g \neq 0}\frac{|a|_p}{\|g\|_p}
\]
which we defined in (\ref{s-ratio}). This quantity exists and is
finite because the set $\{\phi(\bfx\cdot \xi)\}_{\xi\in X}$ is a
linearly independent, finite set of functions. The quantity
$\maxg_{\calg,\,p}$ provides a measure of the linear independence of the
set, albeit one that scales with the norm of $\phi$. Once $\phi$ is
fixed, it depends completely on the geometry of $X$.

For a continuous SBF $\phi$, this is related to the stability of the
interpolation matrix for $\phi$ and $X$.  However, we are only
assuming that $\phi$ is in $L^p$, and thus evaluating $\phi$ on $X$ is
meaningless. Even so, using a smoothed version of $\phi$ allows us to
connect the two concepts.

\subsection{Stability ratios and interpolation matrices}
Let $\kappa \ge 0$ be in $C^k(\RR)$, $k\ge n+2$, and let it satisfy
(\ref{kappa_condits}). Of course, since $\kappa$ is not identically 0,
we also have that there is some open interval on which
$\kappa>0$. Consider the corresponding operator
$\sfk_{\varepsilon,n}=\kappa(\varepsilon \sfl_n)$ and its kernel
$K_{\varepsilon,n}$. To smooth $g(\bfx)=\sum_{\xi\in X}a_\xi \phi(\bfx
\cdot \xi)$, apply $\sfk_{\varepsilon,n}$ to both sides. Doing this
yields
\begin{equation}
\label{g_eps_coef_vec}
g_\varepsilon(\bfx)=\sfk_{\varepsilon,n}g(\bfx) = \sum_{\xi\in X}a_\xi 
\underbrace{\sfk_{\varepsilon,n}
\phi(\bfx \cdot \xi)}_{\displaystyle{\phi_\varepsilon(\bfx \cdot \xi)}}
\end{equation}

We want to relate $\maxg_{\calg,\,p}$ to quantities in a standard SBF
interpolation problem on $X$ involving $\phi_\varepsilon$. The
function $\phi_\varepsilon$ is a spherical harmonic, with nonnegative
Fourier-Legendre coefficients, whose degree depends on the support of
$\kappa$. It is thus a positive definite function on $\sph^n$, but not
an SBF.

The interpolation matrix corresponding to $\phi_\varepsilon$ is 
\[
A_\varepsilon=[\phi_\varepsilon(\eta\cdot \xi)]_{\xi,\eta\in X}.
\]
Later, as a by-product of our analysis, we will establish the
invertibility of $A_\varepsilon$, provided $\varepsilon$ satisfies
certain conditions. When $\varepsilon$ is sufficiently small, one can
also establish it by using a result of Ron and Sun \cite[Theorem
6.4]{Ron-Sun-96-1}: Let $X\subset \sph^n$ be fixed and let $\psi$ be a
positive definite function, but not necessarily an SBF (i.e., some of
coefficients $\hat\psi(\ell)$ may vanish). Then, there is an integer
$j_{X,n}$ such that the interpolation matrix $A_\psi$ will be positive
definite if the set of integers on which $\hat \psi(\ell)>0$ contains
at least $j_{X,n}$ consecutive even integers and $j_{X,n}$ consecutive
odd integers. With our assumptions on $\kappa$ -- in particular, that
$\kappa$ is not identically 0 -- it is clear that for sufficiently
small $\varepsilon$ there are arbitrarily large sets of consecutive
integers for which $\hat \phi_\varepsilon(\ell)>0$. Thus
$A_\varepsilon$ is (strictly) positive definite, and hence invertible,
for all such $\varepsilon$.

Our approach will again be to use the Riesz-Thorin theorem. Let
$y_\varepsilon := g_\varepsilon |_X$, the restriction of
$g_\varepsilon$ to $X$. Using (\ref{g_eps_coef_vec}), we can
interpolate $g_\varepsilon$ on $X$:
\[
y_\varepsilon = A_\varepsilon a\,, \
A_\varepsilon=[\phi_\varepsilon(\eta\cdot \xi)]_{\xi,\eta\in X},
\]
Solving and taking the $\ell^1$ norm, we see that
\[
|a|_1\le \|A_\varepsilon^{-1}\|_1 |y_\varepsilon|_1\,, \
|y_\varepsilon|_1 = \sum_{\xi \in X} |g_\varepsilon(\xi)|.
\]
By our assumptions on $\kappa$ and by (\ref{SUM_function_est}), we
have that $ \ |y_\varepsilon|_1 \le C_{n,\kappa,k} \rho^n
\varepsilon^{-n} \|g\|_{L^1}$. Consequently, for $\phi\in L^1$ we
have that
\[
\maxg_{\calg,\,1} \le C_{\kappa,n,k}\rho^n \varepsilon^{-n}
\|A_\varepsilon^{-1}\|_1\,.
\]
Similarly, working with $p=\infty$ we obtain 
\[
|a|_\infty\le \|A_\varepsilon^{-1}\|_\infty |y_\varepsilon|_\infty\,,
\ |y_\varepsilon|_\infty = \max_{\xi\in X}\{|g_\varepsilon(\xi)|\}\le
\|g\|_\infty.
\]
Recall that $A_\varepsilon^{-1}$ is a self-adjoint matrix, and that
for such matrices the $p=1$ and $p=\infty$ norms are equal:
$\|A_\varepsilon^{-1}\|_\infty = \|A_\varepsilon^{-1}\|_1$. Hence, for
$\phi\in C$ ($p=\infty$), we obtain
\[
\maxg_{\calg,\,\infty} \le \|A_\varepsilon^{-1}\|_1\,.
\]
Applying the Riesz-Thorin theorem to these bounds yields the following:

\begin{prop} \label{stability_ratio_interp_est}
Let $\varepsilon \le \|\calx\|$ and let $\phi\in L^p$. Then, 
\[
\maxg_{\calg,\,p} \le C_{\kappa,n,k}^{1/p}\rho^{n/p} \varepsilon^{-n/p}
\|A_\varepsilon^{-1}\|_1\,.
\]
\end{prop}

\subsection{$\ell^1$ stability estimates for interpolation matrices}

The estimates we need next are for $\|A_\varepsilon^{-1}\|_1$, and the
approach we take to get them will depend on $\phi$ and the behavior of
the $\hat\phi(\ell)$'s.  We will first deal with the Green's function
case, in which $\hat\phi(\ell)$ decays algebraically. After that, we
will deal with the case in which $\phi$ is $C^\infty$, and
$\hat\phi(\ell)$ has very fast decay.

\subsubsection{SBFs that are perturbations of Green's functions}

A straightforward way to estimate the 1-norm of the inverse of a
matrix is to use diagonal dominance techniques, if the matrix is
amenable to them. To that end, split an $n\times n$ matrix $A$ into
its diagonal $D$ and off-diagonal $F$, so $A=D+F$. We then have the
following standard norm estimate, whose proof we omit.

\begin{lemma}
\label{diag_dom_lem}
If $D$ is invertible and $\|D^{-1}F\|_1<1$, then $A$ is invertible and
$\|A^{-1}\|_1<\|D^{-1}\|_1(1- \|D^{-1} F\|_1)^{-1}$.
\end{lemma}

We can apply this to $A_\varepsilon$. The diagonal part is
$D=\phi_\varepsilon(1) I$, and so $\| D^{-1}\|_1=
\phi_\varepsilon(1)^{-1}$ and $\|D^{-1}F\|_1 =
\phi_\varepsilon(1)^{-1}\|F\|_1$. Since the 1-norm of a matrix is the
maximum of the 1-norms of its columns, our condition becomes
\begin{equation}
\label{diag_dom_condit}
\phi_\varepsilon(1)^{-1}\|F\|_1 =
\phi_\varepsilon(1)^{-1}\max_{\eta\in X} 
\sum_{X\ni \xi\neq\eta}| \phi_\varepsilon(\eta\cdot \xi)| <1.
\end{equation}

We now want to deal with a special $\phi_\varepsilon$, which is not
necessarily generated by an SBF $\phi$. Let $\psi$ be a zonal function
in $L^1$, so that
\[
\psi(\xi\cdot\eta)=\sum_{\ell=0}^\infty \hat\psi(\ell)
\frac{\ell+\lambda_n}{\lambda_n\omega_n}\ultra
{\lambda_n} \ell (\xi\cdot \eta ).
\]
We will assume that $1+\hat\psi(\ell)>0$ for all $\ell\ge 0$ and that
$\kappa$ has support in $|t|\in[1,\infty)$. Take $\phi_\varepsilon =
K_{\varepsilon,n}+K_{\varepsilon,n}\ast \psi$, where
$K_{\varepsilon,n}$ is the kernel for the operator $\kappa(\varepsilon
\sfl)$. In addition, define $\psi_\varepsilon=K_{\varepsilon,n}\ast
\psi$.  Since $\hat \phi_{\varepsilon}(\ell) = \kappa(\varepsilon
(\ell+\lambda_n))(1+ \hat\psi(\ell)) \ge 0$, we see that
$\phi_\varepsilon$ is a positive definite spherical function, but not
an SBF. Using (\ref{SUM_kernel_est}) yields
\begin{eqnarray*}
  \sum_{X\ni \xi\neq\eta} |\phi_{\varepsilon}(\eta\cdot \xi)| &\le& 
  \sum_{X\ni \xi\neq\eta} |K_{\varepsilon,n}(\eta\cdot \xi)| + 
  \sum_{X\ni \xi\neq\eta} |\psi_{\varepsilon}(\eta\cdot \xi)| \\
  &\le& C_{n,\kappa,k} q^{-n} + \sum_{\xi\in X}
  |\psi_{\varepsilon} (\eta\cdot \xi)| 
\end{eqnarray*}
Thus, from this and equation (\ref{SUM_function_est}), with
$\kappa(t)=0$, $|t|\le 1$, we have shown that
\begin{equation}
\label{SUM_phi_eps_bnd}
\sum_{X\ni \xi\neq\eta} |\phi_{\varepsilon}(\eta\cdot \xi)| \le 
C_{n,\kappa,k} (q^{-n} + \rho^n\varepsilon^{-n}
E_{L_\varepsilon}(\psi)_1),\ 
L_\varepsilon = \lfloor 1/\varepsilon - \lambda_n \rfloor.
\end{equation}

Thus we have bounded the sum involved in the diagonal dominace
condition (\ref{diag_dom_condit}). Next, we will deal with
$\phi_\varepsilon(1)$. We have the following chain of inequalities:
\begin{eqnarray*}
\phi_\varepsilon(1)&=&K_{\varepsilon,n}(1) +K_{\varepsilon,n}\ast \psi(1) \\
&=&\sum_{\ell=0}^\infty
\kappa(\varepsilon(\ell+\lambda_n))(1+\hat\psi(\ell)) d_\ell^n \\
&\ge&c_0 \sum_{\ell=0}^\infty
\kappa(\varepsilon(\ell+\lambda_n))d_\ell^n = c_0  K_{\varepsilon,n}(1),
\end{eqnarray*}
where $c_0=\min_{\ell\ge 0}(1+\psi(\ell))>0$. (This is true because
$\psi\in L^1$ implies that $\hat\psi(\ell) \to 0$ as $\ell\to
\infty$.) Furthermore, it is easy to see that
\[
K_{\varepsilon,n}(1) = \sum_{\ell=0}^\infty
\kappa(\varepsilon(\ell+\lambda_n)) d_\ell^n \sim
\varepsilon^{-n}\underbrace{  \int_1^\infty\kappa(t)t^{n-1}dt}_{>0}.
\]
Thus, $\phi_\varepsilon(1)\ge C''_{n,\kappa,k} \varepsilon^{-n}$. From
this and (\ref{SUM_phi_eps_bnd}), we arrive at the bound below:
\begin{equation}
\label{diag_dom_phi_eps_bnd}
\|D^{-1}F\|_1\le  C_{n,\kappa,k}\left((\varepsilon/q)^n + \rho^n
  E_{L_\varepsilon}(\psi)_1\right), \ L_\varepsilon = \lfloor
1/\varepsilon - 
\lambda_n \rfloor.
\end{equation}
By choosing $\varepsilon \le q$ sufficiently small, we can make $
C_{n,\kappa,k}\rho^n E_{L_\varepsilon}(\psi)_1$ less than $1/4$, since
$E_{L_\varepsilon}(\psi)_1 \to 0$ as $L_\varepsilon\to \infty$. At
this point, the choice of $\varepsilon$ depends only on $\psi$ and the
mesh ratio $\rho$. If necessary, we may then choose $\varepsilon$
smaller still in order to force the first term on the right to be less
than $1/4$. With this choice of $\varepsilon$, which depends on
$\rho$, $n$, $\kappa$ and $k$, we obtain $\|D^{-1}F\|_1<1/2$. By
Lemma~\ref{diag_dom_lem}, we get the bound on
$\|A_\varepsilon^{-1}\|_1$ below.

\begin{prop} \label{norm_A_1_kernel_est} Suppose that $\kappa$ has
  support in $|t|\in[1,\infty)$. Let $\phi_\varepsilon =
  K_{\varepsilon,n}+K_{\varepsilon,n}\ast\psi$, where $\psi\in L^1$ is
  a zonal function satisfying $1+\psi(\ell)>0$ for $\ell\ge 0$. Then
  there are constants $c$ and $C$, which depend on $\psi$, on $\rho$,
  $n$, $\kappa$ and $k$, such that whenever $\varepsilon \le c q$ we
  have $\|A_\varepsilon^{-1}\|_1 \le C \varepsilon^n $.
\end{prop}

The proof above required conditions on the support of $\kappa$ in
order to deal with the perturbation generated by $\psi$. If $\psi$ is
$0$, then there is no need for such restrictions. Also, the term
involving $\rho$ is gone, and it is no longer involved in determining
$c$ and $C$. We collect these observations below.
\begin{remark}
\label{norm_A_1_remark}
If $\psi=0$, then Propostion \ref{norm_A_1_kernel_est} holds without
restriction on the support of $\kappa$, and neither $c$ nor $C$ depend
on $\rho$.
\end{remark}

We now take an SBF $\phi$ of the form $\phi=G_\beta+G_\beta\ast \psi$,
where $G_\beta$ is the Green's function for $\sfl^{\beta}$ and
$\psi\in L^1$. Our aim is to establish a bound on the stability ratio
for such $\phi$.

\begin{theorem}\label{stability_ratio_algebraic}
  Consider the SBF $\phi=G_\beta+G_\beta\ast \psi$, where $G_\beta$ is
  the Green's function for $\sfl^{\beta}$ and $\psi\in L^1$. Let $X$
  be a set of centers with separation radius $q$ and mesh ratio
  $\rho$. Let $\calg = \calg_{\phi,X}$ be the corresponding SBF
  network. Then there is a constant $C=C(n,\phi,\beta)$ such that
  the stability ratio of $\calg$ satisfies
\begin{equation}
\label{stability_ratio_algebraic_est}
\maxg_{\calg,\,p} \le C\rho^{n/p}q^{n/p'-\beta}
\end{equation}
\end{theorem}

\begin{proof}
  Since we are assuming that $\phi$ is an SBF, the coefficients of the
  $L^1$ function $\psi$ must satisfy $1+\hat\psi (\ell)>0$ for all
  $\ell\ge 0$. Assume $\kappa$ satisfies (\ref{kappa_condits}) and has
  support in $|t|\in[1,\infty)$. The corresponding $\phi_\varepsilon$
  is just $\phi_\varepsilon=\sfk_{\varepsilon,n}\phi =
  \sfk_{\varepsilon,n}(G_\beta + G_\beta\ast\psi)$. By
  Corollary~\ref{kernel_scaling_prop}, we have that
  $\sfk_{\varepsilon,n}G_\beta=\varepsilon^\beta \widetilde
  \sfk_{\varepsilon,n}=\tilde\kappa(\varepsilon \sfl)$, where
  $\tilde\kappa(t)=|t|^{-\beta}\kappa(t)$ satisfies
  (\ref{kappa_condits}). From this, we have that $\phi_\varepsilon =
  \varepsilon^\beta \tilde \phi_\varepsilon$.  If we let $\tilde
  A_\varepsilon$ be the interpolation matrix for $\tilde
  \phi_\varepsilon$, we see that $A_\varepsilon = \varepsilon^\beta
  \tilde A_\varepsilon$. The function $\tilde \phi_\varepsilon$
  satisfies the conditions on the corresponding function in
  Proposition~\ref{norm_A_1_kernel_est}. Thus, by choosing
  $\varepsilon \le cq$, we have
\[
\| A_\varepsilon^{-1}\|_1 = \varepsilon^{-\beta}\| \tilde
A_\varepsilon^{-1}\|_1 \le C\varepsilon^{n-\beta}.
\]
From Proposition~\ref{stability_ratio_interp_est}, we obtain
\[
\maxg_{\calg,\,p}\le C_{\kappa,n,k}^{1/p}\rho^{n/p} \varepsilon^{-n/p}
\|A_\varepsilon^{-1}\|_1\le C' \rho^{n/p} \varepsilon^{n/p'-\beta},
\]
Choosing $\varepsilon$ as large as possible, namely $\varepsilon=cq$,
we have
\[
\maxg_{\calg,\,p}\le C\rho^{n/p} q^{n/p'-\beta},
\]
where the constant $C=C(n,\kappa,k,\phi,p,\beta)$. By taking the
infimum over all $\kappa$, $p$ and $k$, we reduce the dependency of
$C$ to $C=C(n,\phi,\beta)$. This completes the proof.
\end{proof}

\subsubsection{Infinitely differentiable SBFs}
Let $\phi$ be infinitely differentiable SBF. The fast decay of the
Fourier-Legendre coefficient $\hat\phi(\ell)$ requires a different
approach to bounding $\maxg_\calg$ than the one used to obtain
Theorem~\ref{stability_ratio_algebraic}. As before, we let
$A_{\varepsilon}$ be the $N\times N$ interpolation matrix for
$\phi_\varepsilon = \sfk_{\varepsilon,n}\phi$. In addition, we will
let $A$ be the corresponding matrix for $\phi$. By standard matrix
estimates, the norm $\|A_\varepsilon^{-1}\|_1$ satisfies
\[
\|A_\varepsilon^{-1}\|_1\le N^{1/2} \|A_\varepsilon^{-1}\|_2. 
\]
Since $A_\varepsilon$ is a positive definite selfadjoint matrix, the
norm $\|A_\varepsilon^{-1}\|_2$ is equal to the reciprocal of
$\lambda_{\mathrm {min}}(A_\varepsilon)$, the smallest eigenvalue of
$A_\varepsilon$; that is, $\|A_\varepsilon^{-1}\|_2 =
1/\lambda_{\mathrm {min}}(A_\varepsilon)$. We will begin by estimating
this eigenvalue. In preparation for this, we define the quantity
\begin{equation}
\label{hat_phi_min_def}
 \hat\phi_{\mathrm {min}}(L) :=\min_{\,0\le \ell \le L}\hat\phi(\ell)>0.
\end{equation}
where the strict positivity follows from $\phi$ being an SBF.

\begin{prop}\label{min_eval_prop}
  Let $\kappa \ge 0$ be in $C^k(\RR)$, $k\ge n+2$, and let it satisfy
  (\ref{kappa_condits}). In addition, suppose that $\supp(\kappa)
  \subseteq [-2,2]$ and that $\kappa\le 1$. Then, there are constants
  $c=c_{n,\kappa,k}>0$ and $C=C_{n,\kappa,k}>0$ such that for all
  $\varepsilon\le c q$,
\[
\lambda_{\mathrm {min}}(A) \ge \lambda_{\mathrm {min}}(A_\varepsilon) \ge
C\hat\phi_{\mathrm {min}}(L_{\varepsilon/2})\varepsilon^{-n},\
L_{\varepsilon/2} := \lfloor 2/\varepsilon - \lambda_n\rfloor.
\]
\end{prop}

\begin{proof}
Using the Rayleigh-Ritz principle, we thus have
\[
\|A_\varepsilon^{-1}\|_2^{-1}=\lambda_{\mathrm {min}}(A_\varepsilon) =
\min_{a\in \CC^N} a^\ast A_\varepsilon a.
\]
where $A_\varepsilon=[\phi_\varepsilon(\eta\cdot \xi)]_{\xi,\eta\in
  X}$. Because $\phi_\varepsilon$ is a (positive definite) zonal
function, we can use its expansion in spherical harmonics to represent
$\lambda_{\mathrm {min}}(A_\varepsilon)$ via
\begin{equation}
\label{lambda_min}
\lambda_{\mathrm {min}}(A_\varepsilon)= \min_{a\in \CC^N} \left(
  \sum_{\ell=0}^\infty  \sum_{m=1}^{d_\ell}
  \kappa((\ell+\lambda_n)\varepsilon)  \hat\phi(\ell) \bigg|
  \sum_{\xi\in X} Y_{\ell,m}(\xi)a_\xi\bigg|^2\right)
\end{equation}
Since the support of $\kappa$ is $[-2,2]$, the sum above cuts off at
$L_{\varepsilon/2} := \lfloor 2/\varepsilon -
\lambda_n\rfloor$. Consequently, we can bound below $\lambda_{\mathrm
  min}(A_\varepsilon)$ this way:
\[
\lambda_{\mathrm {min}}(A_\varepsilon) \ge \hat\phi_{\mathrm
  min}(L_{\varepsilon/2}) \underbrace{\min_{a\in \CC^N} \left(
    \sum_{\ell=0}^{L_{\varepsilon/2}} \sum_{m=1}^{d_\ell}
    \kappa((\ell+\lambda_n)\varepsilon) \bigg| \sum_{\xi\in X}
    Y_{\ell,m}(\xi)a_\xi\bigg|^2\right)}_{\displaystyle{\lambda_{\mathrm
      min}([K_{\varepsilon,n}(\xi\cdot \eta)])}},
\]
Note that $\lambda_{\mathrm {min}}([K_{\varepsilon,n}(\xi\cdot \eta)])=
\| \, [K_{\varepsilon,n}(\xi\cdot \eta)]^{-1}\|_2^{-1} \le \| \,
[K_{\varepsilon,n}(\xi\cdot \eta)]^{-1}\|_1^{-1}$, because $\|B\|_2\le
\|B\|_1$ for all selfadjoint $B$. The existence of $c$ and $C$ and
their dependencies, along with $\| \, [K_{\varepsilon,n}(\xi\cdot
\eta)]^{-1}\|_1\le C\varepsilon^n$ for $\varepsilon \le cq$, follow
from Proposition~\ref{norm_A_1_kernel_est} and
Remark~\ref{norm_A_1_remark}. Finally, applying the Rayleigh-Ritz
principle, (\ref{lambda_min}), and $0\le \kappa\le 1$, we have that
$\lambda_{\mathrm {min}}(A) \ge \lambda_{\mathrm
  min}(A_\varepsilon)$. This finishes the proof.
\end{proof}

There are two immediate consequences that follow from
Proposition~\ref{min_eval_prop}. The first is a bound on the stability
ratio in this case.

\begin{theorem}\label{stability_ratio_infitely_dif}
  Consider the SBF $\phi$, where $\phi$ is assumend to be infinitely
  differentiable, and let $X$ be a set of centers with separation
  radius $q$ and mesh ratio $\rho$. Let $\calg = \calg_{\phi,X}$ be
  the corresponding SBF network. Then there are positive constants
  $C=C_{n,\kappa,k}$ and $c=c_{n,\kappa,k}$ such that the stability
  ratio of $\calg$ satisfies
\[
\maxg_{\calg,\,p} \le C\rho^{n/p}\frac{q^{n(1/p'-1/2)}}{ \hat\phi_{\mathrm
   {min}}(L_{cq/2})},\ \text{where } L_{cq/2}= \lfloor 2/(cq) -
\lambda_n\rfloor
\]
\end{theorem}

\begin{proof}
  Since $\|A_\varepsilon^{-1} \|_1\le N^{1/2}
  \|A_\varepsilon^{-1}\|_2$, Proposition~\ref{min_eval_prop} implies
  that for $\varepsilon \le cq$,
\[
\|A_\varepsilon^{-1} \|_1\le
C_{n,\kappa,k}\frac{N^{1/2}\varepsilon^n}{\hat\phi_{\mathrm {
    min}}(L_{\varepsilon/2})}
\]
By Proposition~\ref{stability_ratio_interp_est}, we then have that 
\[
\maxg_{\calg,\, p} \le
C_{\kappa,n,k,p}\frac{N^{1/2}\rho^{n/p}\varepsilon^{n/p'}}{\hat\phi_{\mathrm
      min}(L_{\varepsilon/2})}.
\]
Noting that $N\sim q^{-n}$ and choosing $\varepsilon=cq$, which is as
large as possible, we obtain the desired inequality.
\end{proof}

The second consequence is a new stability estimate for interpolation
via a $C^\infty$ SBF $\phi$. Again, let $A$ be the interpolation
matrix for $\phi$ on the set $X$. By Proposition~\ref{min_eval_prop},
$\|A^{-1}\|_2=\lambda_{\mathrm {min}}(A)^{-1} \le
C\varepsilon^n/\hat\phi_{\mathrm {min}}(L_{\varepsilon/2})$. Taking
$\varepsilon=cq$, we obtain a new bound on the norm of $A^{-1}$:
\begin{equation}
\label{new_norm_inverse_bound}
\|A^{-1}\|_2 \le  C\frac{q^n}{\hat\phi_{\mathrm {min}}(L_{cq/2})}.
\end{equation}

\section{Bernstein inequalities and inverse theorems}
\label{bernstein_inverse_theorems}

In this section, we will discuss both direct and inverse theorems for
approximation by SBFs. For an overview of these notions, see
\cite{DeVore-Lorentz-93-1}.

\subsection{Bernstein inequalities}\label{bernstein_inequalities}
Bernstein inequalities are a primary tool in obtaining inverse
theorems. In the introduction, we gave a strategy for obtaining
Bernstein theorems. We have completed the preparation required to state
and prove them. Our first result is for SBFs that are perturbations of
Green's functions.

\begin{theorem}\label{bernstein_greens_fn}
  Consider the SBF $\phi=G_\beta+G_\beta\ast \psi$, where $G_\beta$ is
  the Green's function for $\sfl^{\beta}$ and $\psi\in L^1$. Let $X$
  be a set of centers with separation radius $q$ and mesh ratio
  $\rho$, and let $\calg = \calg_{\phi,X}$ be the corresponding SBF
  network. If $1\le p\le \infty$, $0<\gamma<\beta-n/p'$ and $g\in
  \calg$, then
\begin{equation}
\label{bernstein_algebraic_est}
\|g\|_{H^p_\gamma} \le C q^{-\gamma}\|g\|_p.
\end{equation}
\end{theorem}

\begin{proof}
  Recall that $\|g\|_{H^p_\gamma} \le \|B_Jg\|_{H^p_\gamma} +
  \|(I-\sfb_J)g\|_{H^p_\gamma}$, where $\sfb_J$ is the frame
  reconstruction operator defined in section~\ref{frames}. Of course,
  from (\ref{B_J_approx_id}), this operator is bounded independently
  of $J$ From the polynomial version of the Bernstein inequality in
  (\ref{bernstein_poly_est}), we have that $\|\sfb_J g\|_{H^p_\gamma}
  \le C2^{\gamma J}\|\sfb_J g\|_p \le C2^{\gamma J}\|g\|_p$, which
  implies (\ref{basic_ineq}). Inserting the approximation estimate
  (\ref{final_approx_est}) and the stability-ratio estimate
  (\ref{stability_ratio_algebraic_est}) into (\ref{basic_ineq})
  yields
\begin{eqnarray*}
  \|g\|_{H^p_\gamma} &\le& \left( C2^{\gamma J} + C'
    2^{-(\beta-\gamma-n/p')J}q^{n/p'-\beta}(1+E_{2^{J+j_n}}(\psi)_1 
  \right)\|g\|_p \\
  &\le& q^{-\gamma}\left( C(2^Jq)^{\gamma} +
    C'(2^{-J}q)^{(\beta-\gamma-n/p')}(1+\|\psi\|_1) \right)\|g\|_p
\end{eqnarray*}
The integer $J$ is still a free parameter. Choose it to be
$J=-\log_2(q)$. The Bernstein inequality
(\ref{bernstein_algebraic_est}) then follows on noting that $q\le
\pi$, $\beta-\gamma-n/p'>0$, and $\|\psi\|_1$ is finite and fixed.
\end{proof}

Up to a point, an SBF $\phi\in C^\infty$ is handled in the same way as
one related to a Green's function. In particular, using the argument
above, coupled with the approximation estimate
(\ref{final_approx_est}), with $\beta=\gamma+n$, and the stability estimate in
Theorem~\ref{stability_ratio_infitely_dif}, we obtain
\begin{equation}
\label{bernstein_infty_preliminary}
\|g\|_{H^p_\gamma} \le C L^\gamma \left( 1 +C' 
\rho^n(qL)^{n(1/p'-1/2)}
\frac{L^{-(\beta - \frac{n}{2})}E_L(\sfl_n^\beta \phi)_1}
{ \hat\phi_{\mathrm
    min}(L_{cq/2})} 
\right)\|g\|_p, \ L=2^{J+j_n},
\end{equation}
where $L_{cq/2} =\lfloor 2/cq -\lambda_n\rfloor$. Because $\phi\in C^\infty$, it is in $H^p_\beta$ for all $\beta$. The inequality thus holds for all $\beta>\gamma+n/p'$. The object here is to find a constant $L = \alpha
q^{-1}$, where $\alpha$ is independent of $q$, such that the ratio on
the right above is bounded. The other terms will be controlled easily
in that case. To obtain a simple, applicable condition, we need the following lemma.

\begin{lemma}\label{sequence_condition}
Let $0<\mu(\ell)\le \sigma(\ell)$ be eventually decreasing sequences. Assume that for every $\alpha>0$ there is an integer $m_1=m_1(\alpha,\sigma)\ge 0$ such that  $\ell^{\alpha}\sigma(\ell)\le \sigma(2^{-m_1}\ell)$. If in addition for all $\ell$ sufficiently large  there is an integer $m_2(\alpha,\mu,\sigma)\ge 0$ such that $\sigma(2^{m_2}\ell)\le C_{\mu,\sigma}\mu(\ell)$, then with $m=m_1+m_2$, 
\[
\frac{1}{\mu(L)}\sum_{\ell=2^m L}^\infty \ell^\alpha \sigma(\ell) \le C_{\mu,\sigma}2^{-m} L^{-1}.
\]
\end{lemma}
\begin{proof}
Let $m_1=m_1(\alpha+2,\phi)$. then 
\[
\sum_{\ell=L}^\infty \ell^\alpha \sigma(\ell)\le  \sum_{\ell=L}^\infty \ell^{-2}\ell^{\alpha+2} \sigma(\ell)\le \sigma(2^{-m_1}L)\sum_{\ell=L}^\infty \ell^{-2} \le \frac{\sigma(2^{-m_1}L)}{L}.
\]
Replace $L$ by $2^{m} L$ in the inequality above, so that the sum on the left above is bounded by $(2^{m}L)^{-1}\sigma(2^{m_2}L)\le C_{\mu,\sigma}2^{-m}L^{-1}\mu(L)$. Dividing by $\mu(L)$ yields the desired inequality.
\end{proof}

\begin{lemma}\label{coef_seq_condit}
If there are two sequences $\mu(\ell)$ and $\sigma(\ell)$ that satisfy the conditions of Lemma~\ref{sequence_condition} and in addition satisfy $\mu(\ell)\le \hat\phi(\ell)\le \sigma(\ell)$, then there is an integer $m=m(\beta,\phi,n)$ such that for all $L$ sufficiently large
\begin{equation}
\label{error_bnd_phi_infty_condit}
\frac{E_{2^m L}(\sfl_n^\beta \phi)_1}
{ \hat\phi_{\mathrm {min}}(L)} \le C_{\beta,\phi,n}2^{-m} L^{-1},
\end{equation}
\end{lemma}
\begin{proof}
Because $\phi$ is a $C^\infty$ SBF, the error $E_L(\sfl_n^\beta
\phi)_1$ satisfies 
\[
E_L(\sfl_n^\beta \phi)_1 \le \omega_n
E_L(\sfl_n^\beta \phi)_\infty \le \sum_{\ell=L}^\infty
\frac{(\ell+\lambda_n)^\beta\ultra {\lambda_n} \ell
  (1)}{\lambda_n}\hat\phi(\ell) \le
\frac{2^{\beta+n}}{\Gamma(n)}\sum_{\ell=L}^\infty
\ell^{\beta+n-1} \hat\phi(\ell),
\] 
where we have estimated factors independent of $\phi$ to get the term on the right. Applying Lemma~\ref{sequence_condition} then completes the proof.
\end{proof}

Putting all these results together leads to this theorem.
\begin{theorem}
\label{bernstein_infty_fn}
Let $\phi$ be a $C^\infty$ SBF. If there are two sequences $\mu(\ell)$ and $\sigma(\ell)$ that satisfy the conditions of Lemma~\ref{sequence_condition} and in addition satisfy $\mu(\ell)\le \hat\phi(\ell)\le \sigma(\ell)$, then for every $\gamma>0$ Bernstein's inequality, 
\[
\|g\|_{H^p_\gamma} \le C_{\phi,\gamma,p}
q^{-\gamma} \|g\|_p,
\]
holds for all $g\in \calg_{\phi,X}$, $1\le p\le \infty$. In particular, it holds for the Gaussians, multiquadrics, ultraspherical generating funstions, and the Poisson kernel.
\end{theorem}

\begin{proof}
To get the inequality itself, use Lemma~\ref{bernstein_infty_fn} with $\beta = \gamma+n>\gamma+n/p'$. The statement concerning the list of functions may be established by checking that upper and lower bounds given in section~\ref{SBFs} for each function satisfy the conditions on $\mu(\ell)$ and $\sigma(\ell)$.
\end{proof}

\subsection{Direct theorems} \label{direct_thms}

In \cite[\S 4]{Mhaskar-etal-99-1}, we used a linear process to estimate the distance $\dist_{L^p}(f, \calg_{\phi,X})$, given that $\phi$ is a continuous SBF and $f\in L^p$. In several important cases, including the Gaussian, the process produced a near-best approximant. We will use a similar process here for an SFB of the form $\phi_\beta=G_\beta+G_\beta\ast \psi$, $\psi\in L^1$, again obtaining the corresponding distance estimates.  Such SBFs are at least in $L^1$, but they might not be continuous. Our approach also makes use of recently developed positive-weight quadrature formulas for $\sph^n$, introduced in \cite{Mhaskar-etal-01-1} and further developed in \cite{Narcowich-etal-06-1}.  We remark that a version of Theorem~\ref{phi_beta_favard}, with the conditions on $\phi$ given in terms of sequence spaces involving the $\hat\phi(\ell)$'s, was established in \cite[Theorem~3.1]{Mhaskar-06-1}.

The general framework is this. Let $\phi$ be an SBF, so that the Fourier-Legendre coefficients $\hat\phi(\ell)$ are positive for all $\ell$. Define $\phi^{-1}$ to be the formal expansion
\[
\phi^{-1} \sim \sum_{\ell=0}^\infty 
\frac{\ell+\lambda_n}{\lambda_n\omega_n}\hat\phi(\ell)^{-1}\ultra
{\lambda_n} \ell\,.
\]
This expansion will converge in a distributional sense if the $\hat\phi(\ell)^{-1}$ grow polynomially fast. Otherwise, i.e. for faster growth, the expansion is purely formal. Since we are using it in connection with polynomials of finite degree, this is not a problem.

For every spherical polynomial $S\in \Pi_L$, we can use $\phi^{-1}$ to define an inverse for the convolution operator $S \to \phi\ast S \in \Pi_L$; namely, the expression $\phi^{-1}\ast S$, which is defined by the expansion 
\[
\phi^{-1}\ast S= \sum_{\ell=0}^L \sum_{m=1}^{d_\ell^n}
\frac{\ell+\lambda_n}{\lambda_n\omega_n}\frac{\hat S(\ell,m)}{\hat\phi(\ell)} Y_{\ell,m}
\]
which is just the convolution of $S$ with the polynomial $\sum_{\ell=0}^L\frac{\ell+\lambda_n}{\lambda_n\omega_n}\hat\phi(\ell)^{-1}\ultra {\lambda_n} \ell$. 

Suppose that $S$ is a spherical polynomial for which $\deg S +\lambda_n \le 
2^{J+j_n}$. By Theorem~\ref{bernstein-nikolskii-poly-ineqs}, we have
that $\sfb_JS=S$. In addition, $S= \phi\ast\phi^{-1}\ast S$. Combining these two then yields
\[
S(x) =\sfb_J \phi\ast \phi^{-1}\ast S
= \int_{\sph^n}(\sfb_J\phi)(x\cdot \eta)(\phi^{-1}\ast S)(\eta)d\mu(\eta).
\]
The kernel $\sfb_J\phi(x\cdot \eta)$ is a zonal polynomial with degree less than $2^{J+j_n+1}$. In addition,, $\phi^{-1}\ast S$ is a spherical polynomial of degree $2^{J+j_n-1}$. Thus, the integrand above is a polynomial of degree less than $2^{J+j_n+1} + 2^{J+j_n-1}<2^{J+j_n+2}$.

We will discretize this integral by applying the quadrature formula in \cite[\S 4.2]{Narcowich-etal-06-1}. Let $X$ be a set of centers, with $q$, $h$, $\rho$, and $\calx$ being the separation radius, mesh norm, mesh ratio, and Voronoi  (or similar) decomposition, respectively. Take $L>0$ be an integer. There are positive weights $c_\xi$, $\xi\in X$ and a constant $s_n>0$ (cf. \cite[\S 4.1]{Narcowich-etal-06-1}) such that 
\begin{equation}
\label{quad_formula}
\int_{\sph^n}f(\eta)d\mu(\eta) \doteq \sum_{\xi\in X}c_\xi f(\xi)
\end{equation}
holds exactly for polynomials in $\Pi_L$, provided that $h \le \frac14 s_n^{-1}(L+\lambda_n)^{-1}$. The weights behave like $c_\xi = \calo\left(h^n\right)$, where the constants hidden by ``big'' $\calo$ are dependent only on the dimension $n$. Applying the quadrature formula to the integral representing $S$ yields
\[
S(x) = \sum_{\xi\in X}c_\xi (\sfb_J\phi)(x\cdot \xi)(\phi^{-1}\ast S)(\xi).
\]
Of course we are assuming that $h\sim 2^{-J}$. Let  $\sfq : \Pi_L \to \calg_{\phi,X}$ be given via
\[
\sfq_\calg S(x) := \sum_{\xi\in X}c_\xi \phi(x\cdot \xi)(\phi^{-1}\ast S)(\xi),
\]
and let $g=\sfq_\calg S$, where $\sfq$ is used because of the operator's relationship with quadrature. The difference between $g$ and $S$ is thus
\[
g-S = \sum_{\xi \in X}c_\xi(I-\sfb_J)\phi((\cdot)\cdot \xi) 
(\phi^{-1}\ast S)(\xi) = (I-\sfb_J)g.
\]

We now want to estimate the $H^p_\gamma$ norm of the difference $g-S=(I-\sfb_J)g$ in terms of $\|\phi^{-1}\ast S\|_p$. It is important to note that the norm $\|\phi^{-1}\ast S\|_p$ depends on the degree of $S$ and on $\phi$. We will deal with it later.

The easiest way to estimate  $\|g-S\|_{H^p_\gamma}$ is to employ Theorem~\ref{main_approximation_thm}, where the norm ratios  $\|(I-\sfb_J)g\|_{H^p_\gamma}/|a|_p$ have been estimated. Thus, the task to be accomplished is to relate $|a|_p$ to $ \|\phi^{-1}\ast S\|_p$. To do this, we will again use the Riesz-Thorin theorem.

First of all, we have that $a_\xi$, which is the coefficient of $\phi((\cdot)\cdot \xi)$ in $g$, is given by $a_\xi = c_\xi (\phi^{-1}\ast S)(\xi)$. Thus,  $|a|_\infty = \max_{\xi\in X}c_\xi |(\phi^{-1}\ast S)(\xi)|$. Since $c_\xi = \calo(h^n)$, the bound $|a|_\infty \le C h^n \|\phi^{-1}\ast S\|_\infty$ holds. 

The $p=1$ case  requires more work. Now, $|a|_1 = \sum_{\xi \in X}c_\xi | (\phi^{-1}\ast S)(\xi)|$. 
Since $c_\xi = \calo\left(h^n\right)\le C_n\rho^n q^n\le C''_n \rho^n \min_{\xi\in X}\mu(R_\xi) \le C''_n \rho^n\mu(R_\xi)$, we have
\[
|a|_1\le   C''_n\rho^n \bigg(\sum_{\xi \in X}\mu(R_\xi) |(\phi^{-1}\ast S)(\xi)|\bigg) \le \frac{5 C''_n\rho^n}{4} \|\phi^{-1}\ast S \|_1\,.
\]
The right side above follows on applying the polynomial version of the Marcinkiewicz-Zygmund inequality from \cite[Theorem~4.2]{Narcowich-etal-06-1}, with $\delta = 1/4$, to bound the sum in the middle by $(5/4)\|\phi^{-1}\ast S\|_1$. The Riesz-Thorin theorem then implies
\[
|a|_p \le C_{n,p} \rho^{n/p}h^{n/p'} \|\phi^{-1}\ast S\|_p.
\]
Combining this with the estimate (\ref{final_approx_est}), where $h\sim \varepsilon_J=2^{-(J+j_n)}$ and noting that $g-S=(\sfq_\calg - I)S$, we obtain the following result.

\begin{lemma}\label{direct_est_phi_poly}
Let $\gamma\ge 0$, $1\le p\le \infty$, $\beta>\gamma+n/p'$, $h\sim 2^{-(J+j_n)}$. If $S$ is a spherical polynomial of degree $2^{J+j_n-1}$ or less , then
\begin{equation*}
\|(\sfq_\calg - I)S\|_{H^p_\gamma} \le C_{n,p} \rho^n
h^{\beta-\gamma}\|\phi^{-1}\ast S\|_p
\left\{
\begin{array}{cl}
E_{2^{J+j_n}}(\sfl_n^\beta \phi)_1 & \phi\in H^1_\beta\,,\\
(1+E_{2^{J+j_n}}(\psi)_1 ) & \phi=G_\beta + G_\beta\ast \psi\,.
\end{array}\right.
\end{equation*}
\end{lemma}

\subparagraph{The $\phi_\beta$ case.} We will now focus on the $\phi_\beta$'s. Our immediate concern is estimating $\|\phi_\beta^{-1}\ast S\|_p$.

\begin{lemma}\label{norm_phi_beta_inv_ast_poly}
Let $1\le p\le \infty$, $\beta>0$, $\psi\in L^1$, and $S\in \Pi_L$. If $\phi_\beta=G_\beta+G_\beta\ast \psi$, then there is a constant $C=C_{n,p,\psi}$, which is independent of $\beta$, $L$, and $S$, such that this holds:
\begin{equation}
\label{phi_beta_inv_S}
\|\phi_\beta^{-1}\ast S\|_p \le C_{n,p,\psi} \|S\|_{H^p_\beta}.
\end{equation}
\end{lemma}

\begin{proof}
Note that $\phi_\beta^{-1}\ast S = (\sfl_n^\beta \phi_\beta)^{-1}\ast \sfl_n^\beta S$. The kernel $G_\beta$ is a Green's function for $\sfl_n^\beta$, and so $\sfl_n^\beta \phi_\beta =\delta+\delta\ast \psi=\delta +\psi$, which is to be regarded as a distributional kernel.  Finding $(\sfl_n^\beta \phi_\beta)^{-1}\sfl_n^\beta S$ requires solving $\sfl_n^\beta \phi_\beta\ast T =T+\psi\ast T=\sfl_n^\beta S$ for $T$ in $\Pi_L$, which can be done directly, coefficient by coefficient.  The solution $T$ is of course unique. 

There is another way to look at this equation, in an $L^p$ setting. Suppose that we want to solve $Hf:=f+\psi\ast f=h$ in $L^p$, for $1\le p<\infty$ and in $C$ (for $p=\infty$). The operator norm for $f\to \psi\ast f$ is $\|\psi\|_1$. By Theorem~\ref{B_J_approx_thm}, we have that $\|\psi-\sfb_J\psi\|_1\to 0$ as $J\to \infty$. It follows that the convolution operator with kernel $\psi$ is the norm limit of finite rank operators with convolution kernels, $\sfb_J\psi$. The operator $\psi\ast$ is therefore compact on all $L^p$ and $C$; hence, $Hf=f + \psi\ast f$ has closed range on these spaces. Moreover, a simple coefficient argument shows that $\ker(H)=\{0\}$. The Fredholm Alternative \cite[\S VII.11]{Dunford-Schwartz-58-I} then implies that $\ker(H^\ast)=\{0\}$, so $H^{-1}$ exists and is bounded on all $L^p$ and $C$. Since $\phi_\beta^{-1}\ast S=H^{-1}\sfl_n^\beta S$, we have that
\begin{equation}
\label{phi_beta_inv_H_inv}
\|\phi_\beta^{-1}\ast S\|_p \le \|H^{-1}\|_p \|S\|_{H^p_\beta}.
\end{equation}
We emphasize that $\|H^{-1}\|_p$ is \emph{independent} of $\beta$, $L$, and $S$. It depends only on $p$, $n$, and $\psi$. Consequently, $ C_{n,p,\psi}=\|H^{-1}\|_p$, and (\ref{phi_beta_inv_S}) holds.
\end{proof}

These lemmas lead to the following two direct theorems, the first  for $S\in \Pi_L$ and the second for $f\in H^p_\gamma$.

\begin{theorem}\label{poly_direct_est}
Let $1\le p\le \infty$, $\gamma\ge 0$, and $\beta>\gamma + n/p'$. If $S$ is a spherical polynomial of degree $2^{J+j_n-1}$ or less and if $h=\rho q \sim 2^{-J-j_n}$, then we have for $\phi=\phi_\beta$,
\begin{equation}
\label{phi_beta_poly_direct}
\dist_{H^p_\gamma}(S, \calg_{\phi_\beta,X}) \le C_{n,\beta,\gamma,p,\psi}\rho^n h^{\beta-\gamma} \|S\|_{H^p_\beta}\,, 
\end{equation}

\end{theorem}

\begin{proof} The two lemmas, when applied to $\phi_\beta$, yield
\begin{equation}
\label{norm_est_Q}
\|(\sfq_\calg - I)S\|_{H^p_\gamma} \le C_{n,\beta,\gamma,p,\psi}\rho^n h^{\beta-\gamma}\|S\|_{H^p_\beta}.
\end{equation}
The result follows on observing that $\dist_{H^p_\gamma}(S, \calg_{\phi_\beta,X})\le \|(\sfq_\calg - I)S\|_{H^p_\gamma} $. Note that the dependence of $C$ on the particular frame operator disappears on minimizing the constants involved over all functions $a$.
\end{proof}

\begin{theorem}
\label{phi_beta_favard}
Let $1\le p\le \infty$, $\gamma\ge 0$, and $\beta>\gamma + n/p'$. If $f\in H^p_\beta$, then for $\phi_\beta=G_\beta+G_\beta\ast\psi$, $\psi\in L^1$, 
\[
\dist_{H^p_\gamma}(f, \calg_{\phi,X}) \le C_{\beta,\gamma,n,p,\psi} h^{\beta-\gamma} \rho^n
\|f\|_{H^p_\beta}\,.
\]
\end{theorem}

\begin{proof}
Let $2^{-J-j_n}\sim h$ and choose $S$ to be the polynomial $S=\sfb_J f$; note that $\sfq_\calg S\in \calg_{\phi_\beta,X}$. From these choices and (\ref{norm_est_Q}), it follows that
\begin{align*}
\|f-\sfq_\calg S\|_{H^p_\gamma} &\le \|f-\sfb_Jf\|_{H^p_\gamma}+\|(\sfq_\calg - I)S\|_{H^p_\gamma} \\
&\le \|f-\sfb_Jf\|_{H^p_\gamma} + 
h^{\beta-\gamma}\rho^nC_{\beta, n,p}\|\sfb_J f \|_{H^p_\beta}.
\end{align*}
By Proposition~\ref{H^p_gamma_H^q_beta}, with $p=q$, we have  and 
\[
\|f-\sfb_Jf\|_{H^p_\gamma}\le C_{\beta,\gamma,
  n,a}2^{-(\beta-\gamma) (J+j_n)} E_{2^{J+j_n}}(\sfl_n^\beta f)_p\le  C_{\beta,\gamma,
  n,a}h^{\beta - \gamma} \|f\|_{H^p_\beta}.
\]
From Proposition~\ref{B_J_approx_thm}, we easily see that $\|\sfb_Jf\|_{H^p_\beta} \le C_{\beta,\gamma,n,a}\|f\|_{H^p_\beta}$. Combining all of these inequalities establishes that
\begin{equation}
\label{near_best_approx_H_gamma}
\|f-\sfq_\calg S\|_{H^p_\gamma} \le C_{\beta,\gamma,
  n,a,\psi} \rho^nh^{\beta - \gamma} \|f\|_{H^p_\beta}
\end{equation}
Since $\dist_{H^p_\gamma}(f, \calg_{\phi,X})\le \|f-\sfq_\calg S\|_{H^p_\gamma}$, and since the distance itself doesn't depend on the particular frame function, minimizing over the $a$ yields the result, with the constant independent of $a$.
\end{proof}

\subparagraph{The $C^\infty$ case.} The case in which the SBF $\phi$
is $C^\infty$ was in large part done in
\cite{Mhaskar-etal-99-1}. However, some adjustments need to be made
because the estimates in that paper did not involve $H^p_\gamma$. One
difference is in estimating the norm $\|\phi^{-1}\ast S\|_p$.

\begin{lemma}
\label{norm_phi_infin_inv_ast_poly}
Let $1\le p\le \infty$, $\delta\ge 0$, $L>0$ an integer, and $S\in
\Pi_L$. If $\phi\in H^P_\delta$ is an SBF, then there is a constant
$C=C_n$, depending only on $n$, such that this holds:
\begin{equation}
\label{phi_infin_inv_S}
\|\phi^{-1}\ast S\|_p \le C_n
\frac{L^{n\left|\frac{1}{2}-\frac{1}{p}\right|}}{
  \widehat{\sfl_n^\delta \phi}_{\mathrm {min}}(L)}\|S\|_{H^p_\delta},
\end{equation}
where $\widehat{\sfl_n^\delta \phi}_{\mathrm {min}}(L) = \min_{\,0\le
  \ell \le L}(\ell+\lambda_n)^\delta \hat\phi(\ell)$.
\end{lemma}

\begin{proof}
  We begin by estimating $\|\phi^{-1}\ast S\|_p$. The case in which
  $\phi \in H^p_\delta$ was essentially done in the proof of
  \cite[Theorem~4.1]{Mhaskar-etal-99-1}; the result, which makes use
  of the Nikolskii inequality (\ref{nikolskii_poly_est}), is the
  following. If $S\in \Pi_L$, then the Nikolskii inequality implies
  that
\[
\|\phi^{-1}\ast S\|_p =\|(\sfl_n^\delta \phi)^{-1}\ast \sfl_n^\delta
S\|_p \le C_n L^{n (\frac{1}{2}-\frac{1}{p})_+}\|(\sfl_n^\delta
\phi)^{-1}\ast \sfl_n^\delta S\|_2\,.
\]
At this point, we simply use the $2$-norm estimate done in
\cite[Theorem~4.1]{Mhaskar-etal-99-1} and a second application of
(\ref{nikolskii_poly_est}) to get
\[
\|(\sfl_n^\delta \phi)^{-1}\ast \sfl_n^\delta S\|_2 \le
(\widehat{\sfl_n^\delta \phi}_{\mathrm {min}}(L))^{-1} \|\sfl_n^\delta
S\|_2 \le C_n L^{n (\frac{1}{p}-\frac{1}{2})_+}
(\widehat{\sfl_n^\delta \phi}_{\mathrm {min}}(L))^{-1} \|\sfl_n^\delta
S\|_p\,.
\]
Putting the two inequalities together completes the proof.
\end{proof}

Let $\phi\in C^\infty$. We can now estimate the $H^p_\gamma$ distance
of $S\in \Pi_L$ to $\calg_{\phi,X}$, in terms of $\|S\|_{H^p_\delta}$,
where $\delta>\gamma+n/p'$. In Lemma~\ref{direct_est_phi_poly}, let
$\beta=\delta+n/2$. Apply Lemma~ \ref{norm_phi_infin_inv_ast_poly},
noting that $L\le 2^{J+j_n-1}\le h^{-1}$ implies
$L^{n\left|\frac{1}{2}-\frac{1}{p}\right|}\le L^{n/2} \le
h^{-n/2}$ to get this:
\begin{equation}
\label{dist_poly_norm_est_Q}
\dist_{H^p_\gamma}(S, \calg_{\phi,X}) \le \|(\sfq_\calg - I)S\|_{H^p_\gamma}
\le C_{n,p} \rho^n
h^{\delta-\gamma}\frac{E_{2^{J+j_n}}(\sfl_n^{\delta+n/2}
  \phi)_1}{\widehat{ \sfl_n^\delta \phi}_{\mathrm {min}}(L)}
\|S\|_{H^p_\delta}\,.
\end{equation}

\begin{theorem}\label{phi_H^p_alpha__poly_favard}
  Let $1\le p\le \infty$, $\gamma\ge 0$, $\delta>\gamma + n/p'$, and
  $\phi\in C^\infty$. If there is an integer $m=m(\delta,\phi)>0$ such
  that
\begin{equation}
\label{phi_H^p_gamma_error_ratio}
\sup_{\ell>0}\frac{E_{2^m \ell}(\sfl_n^{\delta+ n/2} \phi)_1}{ 
  \widehat{\sfl_n^\delta \phi}_{\mathrm {min}}(\ell)}  \le C_{m,n,\delta,\phi}
\end{equation}
holds, and if $S\in \Pi_L$, with $L \le 2^{J+j_n-1-m}$ and $h \sim
2^{-J -j_n}$, then
\begin{equation}
\label{phi_H^p_delta_poly_direct}
\dist_{H^p_\gamma}(S, \calg_{\phi,X}) \le C_{m,n,p,\delta,\gamma} 
h^{\delta-\gamma}\rho^n  \|S\|_{H^p_\delta}\,.
\end{equation}
In addition, for $f\in H^p_\gamma$, we have that
\begin{equation}
\label{phi_H^p_beta_H^p_gamma_direct}
\dist_{H^p_\gamma}(f, \calg_{\phi,X}) \le
C_{m,n,p,\gamma,\delta,\phi} h^{\delta-\gamma} \rho^n \|f\|_{H^p_\delta}.
\end{equation}
Finally, these estimates hold for Gaussians, multiquadrics,
ultrasherical generating functions and Poisson kernels.
\end{theorem}

\begin{proof}
If (\ref{phi_H^p_gamma_error_ratio}) holds, then, since $\widehat{\sfl_n^\delta \phi}_{\mathrm {min}}(L) \ge \widehat{\sfl_n^\delta \phi}_{\mathrm {min}}(2^{J+j_n-m})$, it follows that
\[
\frac{E_{2^{J+j_n}}(\sfl_n^{\delta+n/2} \phi)_1}{\widehat{\sfl_n^\delta \phi}_{\mathrm {min}}(L)} \le \frac{E_{2^{J+j_n}}(\sfl_n^{\delta+n/2} \phi)_1}{\widehat{\sfl_n^\delta \phi}_{\mathrm {min}}(2^{J+j_n-m})}\le C_{m,n,\delta,\phi},
\]
and (\ref{phi_H^p_delta_poly_direct}) follows form this and (\ref{dist_poly_norm_est_Q}). One can establish the $H^p_\gamma$ distance estimate (\ref{phi_H^p_beta_H^p_gamma_direct}) below using a proof virtually identical to that for Theorem~\ref{phi_beta_favard}. Essentially the same argument used in section~\ref{bernstein_inequalities} can be used  here to show that Gaussians, multiquadrics, etc.\ satisfy (\ref{phi_H^p_gamma_error_ratio}), and so the estimates hold for them, too.
\end{proof}

\subsection{Besov spaces.}\label{besov_spaces}
In this section, we review the definitions and basic facts regarding
Besov spaces on $\sph^n$. These spaces, which will interpolate between
$L^p(\sph^n)$ and $H^p_\gamma$, are defined in
\cite{Triebel-86-1}. Other, equivalent definitions of Besov spaces on
$\sph^n$ are given \cite{Narcowich-etal-06-2}. Below, we will make use
of a general construction general construction found in
\cite[Chapters~6]{DeVore-Lorentz-93-1} to characterize these spaces in
terms of spaces of SBF networks, $\calg_{\phi,X}$.

There are two ingredients. First, we need to introduce
certain sequence spaces. If $r>0$ and $0<\tau\le \infty$, we define
for a sequence ${\bf a}=\{a_n\}_{n=0}^\infty$ of real numbers,
\begin{equation}\label{seqbdef}
 \|{\bf a}\|_{\tau,r} := \left\{\begin{array}{ll}\displaystyle
 \left\{\sum_{n=0}^\infty 2^{nr\tau}|a_n|^\tau\right\}^{1/\tau},
 & \mbox{if } 0<\tau<\infty,\\
 \displaystyle\sup_{n\ge 0} 2^{nr}|a_n|, & \mbox{if }
 \tau=\infty.\end{array}\right.
\end{equation}
The space of sequences ${\bf a}$ for which $\|{\bf
  a}\|_{\tau,r}<\infty$ will be denoted by $\seqb_{\tau,r}$.

The other ingredient in the definition of Besov spaces is a $K$--functional \cite[Chapter~6]{DeVore-Lorentz-93-1}. For $\delta, \gamma>0$, $1\le p\le\infty$ and $f\in L^p$, the $K$--functional for $L^p$ and $H^p_\gamma$ is given by
\begin{equation}\label{kfuncdef}
\K_\gamma(p,f,\delta):=\inf_{g\in H^p_\gamma}\{\|f-g\|_p +\delta^\gamma (\|g\|_p+\| g\|_{H^p_\gamma})\}.
\end{equation}
If $r>0$, $0<\tau\le \infty$, $r<\gamma$, we define the class of all $f\in L^p$ for which 
\begin{equation}\label{besovnormdef}
\|f\|_{r,\gamma,\tau,p}:=\|f\|_p + \|\{\K_\gamma(p,f,2^{-n})\}_{n=0}^\infty\|_{\tau,r} <\infty.
\end{equation}
to be the Besov space $B^r_{\tau,p}$. As we shall see, other than the requirement $r<\gamma$, the $\gamma$ dependence will disappear from the characterization of the space, so it isn't necessary to keep it in designating the space.

An important problem in approximation theory is to characterize Besov spaces using degrees of approximation of functions. We recall the results \cite[Theorems~7.5.1 and 7.9.1]{DeVore-Lorentz-93-1}, as it applies in the context of the present paper.

\begin{prop}\label{besovprop}
Let $1\le p \le \infty$, $\gamma>0$, and let $\{V_j\}_{j=0}^\infty$, with $V_0=\{0\}$, be a nested sequence of finite dimensional linear subspaces of $L^p$, $p<\infty$ or $C$, $p=\infty$
Suppose that for $j=1,2,\cdots$, one has both the Favard (Jackson) estimate
\begin{equation}\label{genfavard}
\dist_{L^p}(f,V_j) \le C\,2^{-j\gamma}(\|f\|_p+\|f\|_{H^p_\gamma}),
\end{equation}
for all $f \in H^p_\gamma$, and the
Bernstein inequality
\begin{equation}\label{genbernineq}
\|g\|_{H^p_\gamma} \le C 2^{j\gamma}\|g\|_p, \qquad g\in V_j.
\end{equation}
Then for $0<r<\gamma$, $0<\tau\le\infty$, $f\in B^r_{\tau,p}$ if and only if $\{\dist_{L^p}(f,V_j)\}_{j=0}^\infty\in \seqb_{\tau,r}$. 
\end{prop}

\begin{proof}
  This is just \cite[Theorem~7.5.1]{DeVore-Lorentz-93-1}, with the
  sequence of spaces satisfying all requirements in listed in
  \cite[(5.2), p.~216]{DeVore-Lorentz-93-1}, except possibly
  density. This requirement is in fact satisfied if the Favard
  inequality (\ref{genfavard}) is satisfied. To see this, note that
  $H^p_\gamma$ contains all of the spherical polynomials, which form a
  dense set in $L^p$, $1\le p<\infty$ and in $C$. The Favard
  inequality (\ref{genfavard}) then implies that the $\cup_j V_j$ is
  dense in $H^p_\gamma$ and therefore in $L^p$, $1\le p<\infty$, or in
  $C$.
\end{proof}

In the important case when $V_j=\Pi_{2^j}$,
Proposition~\ref{bernstein-nikolskii-poly-ineqs} gives the Bernstein
estimate, while Proposition~\ref{H^p_gamma_H^q_beta} provides the
Favard estimate. In addition, since the criterion that
$\{\dist_{L^p}(f,\Pi_{2^J})\}_{n=0}^\infty\in \seqb_{\tau,r}$ does
\emph{not} depend upon $\gamma$, it follows that the Besov spaces
$B^r_{\tau,p}$ are independent of the different choices of $\gamma>r$
in their definition. This is why we don't need to include the
parameter $\gamma$ to index these spaces.

\begin{remark}
  The polynomial characterization of $B^r_{\tau,p}$ is precisely the
  one given in \cite[Proposition~5.3]{Narcowich-etal-06-2}, so that
  the ``needlet'' definition
  \cite[Definition~5.1]{Narcowich-etal-06-2} is equivalent to the one
  above. (See also \cite{Mhaskar-Prestine-05-1}.) The needlet
  definition is itself known to be equivalent (cf.\
  \cite{Narcowich-etal-06-2}) to that given in \cite{Triebel-86-1}. It
  follows that all three are equivalent.
\end{remark}

Using the proposition above, one can also characterize Besov spaces
using a variety of spherical basis functions. To do this, we must
first have an appropriate nested sequence of sets of centers. By
Proposition~\ref{nesting_Xs}, we can find a nested sequence
$\{X_j\}_{j=0}^\infty \in \calf_\rho$, $\rho\ge 2$, each $X_j$ having
mesh norm $h_j:=h_{X_j}$ satisfying $\frac14 h_j < h_{j+1}\le \frac12
h_j\le \frac{1}{2^j}h_0$. If $\phi\in L^p$ is an SBF, then define the
$V_j$'s to be
\begin{equation}
\label{nested_networks}
V_j:=\calg_{\phi,X_j}, \ j=1, 2, \ldots,\ \mbox{and} \ V_0=\{0\}.
\end{equation}
These spaces have finite dimension equal to the cardinality of $X_j$
and by virtue of the $X_j$'s being nested, are themselves nested. At
issue then are the Favard and Bernstein inequalities. Since any $\phi$
that satisfies both will provide us with a Besov space via
Proposition~\ref{besovprop}, we have the following result.

\begin{cor} \label{besov_cor_phi} Let $1\le p\le \infty$,
  $\phi_\beta=G_\beta+G_\beta\ast \psi $, where $\psi\in L^1$ and
  $0<\beta$. Fix $0<\gamma<\beta-n/p'$ and suppose that
  $V_j=\calg_{\phi_\beta,X_j}$, with $X_j$ as in
  (\ref{nested_networks}). For all $0<r<\gamma$ and all
  $0<\tau\le\infty$, we have that $f\in B^r_{\tau,p}$ if and only if
  $\{\dist_{L^p}(f,V_j)\}_{j=0}^\infty\in \seqb_{\tau,r}$. The same
  conclusion holds true, with any $\gamma>0$, for all $\phi$ that
  simultaneously satisfy (\ref{error_bnd_phi_infty_condit}) and
  (\ref{phi_H^p_gamma_error_ratio}), including the Gaussians,
  multiquadrics, etc.
\end{cor}

\begin{proof}
When $V_j=\calg_{\phi_\beta,X_j}$, the result follows immediately from the Bernstein inequality in Theorem~\ref{bernstein_greens_fn} and the Favard inequality in Theorem~\ref{phi_beta_favard}. 
If $\phi$ satisfies both (\ref{error_bnd_phi_infty_condit}) and (\ref{phi_H^p_gamma_error_ratio}), then it also satisfies both the Bernstein inequality in Theorem~\ref{bernstein_infty_fn} and the Favard inequality in Theorem~\ref{phi_H^p_alpha__poly_favard}. As before, with the same set of $V_j$'s, the same conclusion holds.
\end{proof}

\subsection{Inverse theorems}\label{inverse_theorems}
Inverse theorems give indications of rates of approximation being
best, or nearly best, possible. We now establish inverse theorems for
the approximation rates in the previous section and in \cite{Mhaskar-etal-99-1}. These involve Bessel-potential Sobolev spaces, and in addition Besov spaces.

\begin{theorem}\label{inverse_thm}
  Let $1\le p\le \infty$ and let $\phi $ as in
  Theorem~\ref{bernstein_greens_fn} or
  Proposition~\ref{bernstein_infty_fn}. If for $f\in L^p$, $1\le
  p<\infty$, or $f\in C(\sph^n)$, $p=\infty$, there are constants
  $0<\mu\le \gamma$, $t \in \RR$, and $c_f>0$ such that
\begin{equation}\label{dist_f_bnd}
\dist_{L^p(\sph^n)}(f,\calg_{\phi,X}) \le c_f\frac{ h_X^\mu}{\log_2^t (h_X^{-1})}
\end{equation}
holds for all $X\in \calf_\rho$, then, for every $0\le \nu <\mu$,
$f\in H^p_\nu(\sph^n)$. If  (\ref{dist_f_bnd}) holds for $\nu=\mu$ and some $t>1$, then $f\in H^p_\mu(\sph^n)$. Moreover, if in addition $\phi$ satisfies the conditions in Corollary~\ref{besov_cor_phi}, then for any $\tau>t^{-1}>0$ and $0<r\le \mu$, the function $f$ is in the Besov space $B^r_{\tau,p}$.
\end{theorem}

\begin{proof}
Let the $V_j$'s be as in (\ref{nested_networks}), and set  $f_j:= \argmin \left(\dist_{L^p(\sph^n)}(f,V_j)\right)$, which always exits because $V_j$ is finite dimensional. Since the $V_j$'s are nested, we have that $f_j\in V_k$ for all $k\ge j$.  We want to show that $f_j$ is a Cauchy sequence in $H^p_\nu$. From
  the Bernstein estimate in Theorem~\ref{bernstein_greens_fn} -- or
  Proposition~\ref{bernstein_infty_fn} -- and the inequality
  $h_{j+1}/q_{j+1}\le \rho$, we have
\[
\|f_{j+1}-f_j\|_{H^p_\nu} \le C \rho^\nu
h_{j+1}^{-\nu}\|f_{j+1}-f_j\|_p\le C\rho^\nu h_{j+1}^{-\nu}\big(\|f_{j+1}-f\|_p +
  \|f-f_j\|_p\big).
\]
And by (\ref{dist_f_bnd}), we also have 
\begin{align*}
  \|f_{j+1}-f_j\|_{H^p_\nu}   &\le Cc_f \rho^\nu h_{j+1}^{-\nu} (h_{j+1}^\mu\log_2^{-t}(h_{j+1}) +
  h_j^\mu\log_2^{-t}(h_j))\\
  &\le  Cc_f \rho^\nu h_0 2^{-(\mu - \nu)(j+1)}\left( (h_0+j+1)^{-t}
    + 2^\mu (h_0+j)^{-t}\right)\\
  &\le  C'c_f 2^{-(\mu - \nu)j} j^{-t}
\end{align*}
where $C'$ is independent of $j$. Take $k>j$. Using the previous
inequality and a standard telescoping-series argument, we arrive at
this:
\[
\|f_j-f_k\|_{H^p_\nu} \le C''(\sum_{m=j}^k 2^{-(\mu - \nu)m} m^{-t}).
\]
Letting $j,k\to\infty$, we see $\|f_j-f_k\|_{H^p_\nu} \to 0$ when
$\mu>\nu$ and $\tau \in \RR$ or when $\mu=\nu$ and $t >1$ .  Thus,
$f_j$ is a Cauchy sequence in $H^p_\nu $ and is therefore convergent
to $\tilde f\in H^p_\nu$. Moreover, by (\ref{dist_f_bnd}) with
$X=X_j$, we see that $f_j \to f$ in $L^p$, so $\tilde f = f$ almost
everywhere. Hence, we have $f\in H^p_\nu$. The statement concerning
Besov spaces follows from two things: the observation that
$a_j:=\dist_{L^p}(f,V_j)\le c_f 2^{-\mu j}j^{-t}$, so $\| \bfa
\|_{\tau,r}<\infty$ whenever $0<r \le \mu$ and $\tau t>1$, and
Corollary~\ref{besov_cor_phi}.
\end{proof}

For the case $\nu=\mu$, $0<t\le 1$, the inverse theorem fails for
Bessel-potential Sovolev spaces, but still remains valid for Besov
spaces with $\tau>t^{-1}$.

\section{Concluding Remarks}

There are connections between this paper and \cite{Mhaskar-etal-99-1,
  Mhaskar-06-1}. In these papers quasi interpolatory SBF networks were
obtained yielding near best approximants for functions in Sobolev
classes. The associated quasi-interpolation operators were constructed
in the Fourier domain. The paper \cite{Mhaskar-etal-99-1} focused on
sequences corresponding to the $c^{\infty}$ case treated within this
paper. The paper \cite{LeGia-Mhaskar-06-1} dealt with sequences
connected to the ``perturbations of Green's functions'' case. For
example, let $\psi$ be a perturbation of a Green's function as
described in this paper. If the Fourier coefficents of $\psi$ satisfy
the ``difference condition'' as stated in , \cite{Mhaskar-06-1} then
it is in $L^1$. The examples given in Section 3 satisfy both kinds of
conditions.

In \cite{Mhaskar-etal-99-1, Mhaskar-06-1}, the quasi-interpolatory SBF
networks were shown to give best results in the sense of n-widths. In
this paper, using the frame approach, we have shown the
quasi-interpolatory networks are also optimal for approximation of
individual functions. Also note that in \cite{Mhaskar-06-1},
Marcinkiewicz-Zygmund measures generalizing the measure that
associates $\mu_q(R_\xi)$ with each $\xi$ were introduced. These
measures were used to derive \cite[Prop. 4.1 \& (4.15)]{Mhaskar-06-1},
which have overlap with the current Prop. 4..4, Lemma 4.7 and estimate
(5.4).

In \cite{LeGia-Mhaskar-06-1}, the quasi-interpolation polynomial
operators were further utilized to show that, in the presence of
certain singularities, they exhibited better approximation properties
than traditional methods. Also \cite[Prop. 4.3]{LeGia-Mhaskar-06-1} is
related to Proposition 4.10 given here.  Finally there is material
closely connected to Theorem 4.1 appearing in [6, \cite[Proposition
4.1]{Mhaskar-05-2}.  Another version of the operator $B_J$ was
introduced in \cite{Mhaskar-etal-00-1}: $\sigma_J(f) = \sum_{l =
  0}^{2^J} h(l/2^J)P_l(f),$ where $h:[0, \infty) \rightarrow [0,
\infty)$ is a function in $C^k$, equal to 1 on [0, 1/2] and 0 on $[1,
\infty)$. An early form of Theorem 4.1 was Theorem 3.4 of
\cite{Mhaskar-etal-00-1}.  Frames, based on the $\sigma_J(f)$ operator
can be constructed as in \cite{Mhaskar-05-1, Mhaskar-Prestine-05-1}
using $h(t) - h(2t)$ in place of ${\cal \kappa}$ used in the
construction given here.

Finally we mention that the idea of using minimal separation for
converse theorems and Bernstein inequalities goes back to
\cite{Schaback-Wendland-02-1}, see also \cite{Mhaskar-05-1}.  Also,
for the neural network community, we note that the number of neurons
is not used as a measure of complexity, but rather the minimal
separation of the nodes.

%\bibliography{rbf}

\end{document}